%% file: Articolo_r1_v12.tex
\RequirePackage{fix-cm}
\documentclass[smallextended]{svjour3}       % onecolumn (second format)
\smartqed  % flush right qed marks, e.g. at end of proof
\usepackage[english]{babel}
\usepackage[T1]{fontenc}
\usepackage{graphicx} 
\usepackage{soul}		% Hyphenation for letterspacing, underlining, and more.

\usepackage{psfrag}
\usepackage{subcaption}
\usepackage{color}
\usepackage{amsmath}
\usepackage{amsfonts}
\usepackage{amssymb}
\usepackage{cases}
\usepackage{cancel}
\usepackage{caption}
\usepackage{mathcomp}
\usepackage{multirow}
\usepackage{array}		     	%
\usepackage{hhline}
\usepackage{float}
\usepackage{booktabs}
\usepackage[amssymb]{SIunits}
\usepackage[applemac]{inputenc}
\usepackage[ruled,vlined,linesnumbered]{algorithm2e}
\usepackage{a4wide}
\usepackage{lineno}			%\linenumbers
\usepackage[normalem]{ulem} 	% underline
\usepackage{color}
\usepackage{xcolor}
\usepackage{hyperref}

\usepackage{natbib}
\setcitestyle{authoryear}
\graphicspath{{./figure/}}

\usepackage{xcolor}

\input{macro_ulisse_Springer}

\begin{document}

\title{{Mixed variational formulations for structural topology
  optimization based on  the phase-field approach}}
 
\author{Michele Marino \and Ferdinando Auricchio \and Alessandro Reali \and Elisabetta Rocca \and Ulisse Stefanelli}

%\authorrunning{Short form of author list} % if too long for running head

\institute{M. Marino\at
             Dipartimento di Ingegneria Civile e Ingegneria Informatica, Universit\`a degli Studi di Roma Tor Vergata, via del Politecnico 1, 00133 Roma, Italy 
%              \email{fauthor@example.com}           %  \\
%             \emph{Present address:} of F. Author  %  if needed
           \and
           F.Auricchio (Corresponding), A.Reali \at
            Dipartimento di Ingegneria Civile e Architettura, Universit\`a degli Studi di Pavia, via Ferrata 3, I-27100 Pavia, Italy\\
          Istituto di Matematica Applicata e Tecnologie Informatiche {\it E. Magenes}, via Ferrata 1, I-27100 Pavia, Italy\\
             \email{auricchio@unipv.it}    
              \and
           E.Rocca \at
           Dipartimento di Matematica, Universit\`a degli Studi di Pavia, via Ferrata 3, I-27100 Pavia, Italy
            \and
          U.Stefanelli \at
          Faculty of Mathematics, University of Vienna, Oskar-Morgenstern-Platz 1, A-1090 Vienna, Austria, \\
          Vienna Research Platform on Accelerating Photoreaction Discovery, University of Vienna, W\"ahringerstra\ss e 17, 1090 Wien, Austria\\
          Istituto di Matematica Applicata e Tecnologie Informatiche {\it E. Magenes}, via Ferrata 1, I-27100 Pavia, Italy
}

\date{Received: date / Accepted: date}
% The correct dates will be entered by the editor

\maketitle

\begin{abstract}
 We propose  a  variational principle combining a phase-field functional for structural topology optimization  with a mixed (three-field) Hu-Washizu functional, then including directly in the formulation equilibrium, constitutive, and compatibility equations. 
The resulting mixed variational functional is then specialized  to
derive a classical topology optimization formulation (where the amount
of material to be distributed is an \emph{a priori} assigned quantity  acting as a global constraint for the problem)  as well as a novel topology optimization formulation (where the amount of material to be distributed is minimized, hence with no pre-imposed constraint for the problem). 

Both formulations are numerically solved  by  implementing a mixed finite element scheme, with the second approach avoiding the introduction of a global constraint, hence respecting the convenient local nature of the finite element discretization. Furthermore, within the proposed approach 
it is possible to obtain guidelines for settings proper values of phase-field-related simulation parameters and, 
thanks to the combined phase-field and Hu-Washizu rationale, a monolithic algorithm solution scheme can be easily adopted. 

An insightful and extensive numerical investigation results in a detailed convergence study and a discussion on the obtained final designs. 
The numerical results  clearly highlight differences between the two formulations as well as 
advantages related to the monolithic solution strategy;
numerical investigations address both two-dimensional and three-dimensional applications.

% \PACS{PACS code1 \and PACS code2 \and more}
% \subclass{MSC code1 \and MSC code2 \and more}
\end{abstract}

\input{opt_intro_r1_v12}

		 \input{opt_cont_r1_v12}

		 \input{opt_discretization_r1_v12}
		 		 \input{opt_results_r1_v12}

\clearpage 
\input{opt_conclusion_r1_v12}

\section*{Conflict of interest}
The authors declare that they have no conflict of interest.

\section*{Replication of Results}
The paper already contains all the data necessary to properly replicate all the presented results.

%\newpage

\bibliography{bibliography}

\clearpage 

\appendix

\input{opt_appendix_r1_v12}

%\begin{acknowledgements}
%If you'd like to thank anyone, place your comments here
%and remove the percent signs.
%\end{acknowledgements}

% BibTeX users please use one of
%\bibliographystyle{spbasic}      % basic style, author-year citations
%\bibliographystyle{spmpsci}      % mathematics and physical sciences
%\bibliographystyle{spphys}       % APS-like style for physics
%\bibliography{}   % name your BibTeX data base

%% Non-BibTeX users please use
%\begin{thebibliography}{}
%%
%% and use \bibitem to create references. Consult the Instructions
%% for authors for reference list style.
%%
%\bibitem{RefJ}
%% Format for Journal Reference
%Author, Article title, Journal, Volume, page numbers (year)
%% Format for books
%\bibitem{RefB}
%Author, Book title, page numbers. Publisher, place (year)
%% etc
%\end{thebibliography}

\end{document}

%% file: macro_ulisse_Springer.tex
\newcommand{\blue}     {\color{blue}}

\definecolor{lightgray}{gray}{0.5}

\definecolor{mydgreen}{rgb}{0,0.4,0.3}

\newcommand{\htext }[2]    {\hspace{#1} \text{#2} \hspace{#1}}

\newcommand{\UUU}{\color{blue}}
\newcommand{\OOO}{\color{black}}
\newcommand{\bbb}   {\boldsymbol{b}}

\newcommand{\bbn}   {\boldsymbol{n}}

\newcommand{\bbs}   {\boldsymbol{s}}
\newcommand{\bbt}   {\boldsymbol{t}}
\newcommand{\bbu}   {\boldsymbol{u}}

\newcommand{\bbx}   {\boldsymbol{\mathit{x}}}

\newcommand{\bbN}   {\boldsymbol{N}}

\newcommand{\bbR}   {\boldsymbol{R}}

\newcommand{\bbT}   {\boldsymbol{T}}

\newcommand{\bbX}   {\boldsymbol{X}}

\newcommand{\fourC}      {\mathbb{C}}

\newcommand{\real}       {{\mathcal{R}}}                      
 
\newcommand{\bbzero}  {\boldsymbol{0     }}

\newcommand{\ffxi}  {\boldsymbol{\xi     }}
\newcommand{\ffe}   {\boldsymbol{\varepsilon}}

\newcommand{\ffs}   {\boldsymbol{\sigma  }}

\newcommand{\cC}    {\mathcal{C}}
\newcommand{\cE}    {\mathcal{E}}

\newcommand{\cJ}    {\mathcal{J}}
\newcommand{\cL}    {\mathcal{L}}

\newcommand{\cP}    {\mathcal{P}}

\newcommand{\diver       } {\text{div }}

\newcommand{\grads}      	  {\nabla^s}
\newcommand{\gradsu}      {\nabla^s \bbu}

\newcommand{\bbuhatk}[1]  {\hat{\bbu}_{#1}}

\newcommand{\kphi}      		{\kappa_{\phi}}
\newcommand{\gammaphi}      {\gamma_{\phi}}

%% file: opt_intro_r1_v12.tex
\section{Introduction}

Having assigned a design region, a load distribution, and suitable boundary conditions, the goal of structural 
\textit{topology optimization} is to identify an optimal distribution of material within the design region. Classically, optimality is reached when the obtained material distribution minimizes a measure of structural compliance and satisfies mechanical equilibrium. The minimization of structural compliance is in general expressed in terms of the work done by the assigned (assumed to be constant) load distribution by the corresponding  displacements.

Several methods have been developed in order to solve topology optimization problems.  Originally,  a discrete formulation  was introduced where areas of dense material and voids are alternated without any transition region \citep{bendsoe_83}. 
Known as \textit{0-1 topology optimization} problem, this first approach also leads to many difficulties both from  an analytical and a numerical point of view \citep{sigmund_98}.

Possible alternative approaches are based on homogenization methods, where the macroscopic material properties are obtained from microscopic porous material characteristics \citep{allaire_2004,suzuki_91}, or on the \textit{Solid Isotropic Material Penalization} (SIMP) method \citep{zhou_91}, which consists of penalizing the density region, different from the void or bulk material, by choosing a suitable interpolation scheme for material properties at the macroscopic scale
\citep{bendsoe_99,bendsoe_83}.
Another approach used by many authors  \citep[cf., e.g.,][]{burger_2003,osher_santosa_01} is the level set method, that allows for topology changes but presents difficulties in creating holes. 

An alternative to  SIMP and  level set  is based on the phase-field method, for the first time introduced by~\cite{bourdin_2003} and successively employed by~\cite{burger_2006} for stress constrained problems as well as  by~\cite{takezawa_2010} for minimum compliance and eigenfrequency maximization problems. More recently,~\cite{penzler_2012} have solved nonlinear elastic problems by means of the phase-field approach, while~\cite{dede_isogeometric_2012} have applied this method in the context of isogeometric topology optimization. 

 The  phase-field approach  delivers an   efficient method to solve the free boundary problem stemming from the fact that the boundary of the region filled by the material is unknown. Phase-field topology optimization penalizes an approximation of the interface perimeter in such a way that, by choosing a very small positive penalty term, one can obtain a sharp interface region separating solid materials and voids \citep{blank_multi-material_2014}.
The phase-field approach to structural optimization problems has been recently used by different authors  \citep[cf., e.g.,][]{auricchio_19, carraturo_19, blank_14}, the main advantage being the fact that it allows to handle topology changes as well as nucleation of new holes.
%
%In fact, by inspecting the behavior of the perimeter as a function of the material parameters, one can actually assess the quality of the obtained solutions. 
%Indeed, sudden variations of the perimeter correspond to topology changes in the solution and tend to occur at   critical values of the parameters. The identification of such critical regimes is paramount to the robustness of the topology optimization process, for topological changes corresponding to small parameter variations are clearly undesired. 
%
Moreover, especially for three-dimensional applications, another advantage is that the phase-field approach produces regular patterns in the design domain that require little post-processing effort to interpret results \citep{Deaton2014}. On the other hand, density-based topology optimization approaches (e.g., based on SIMP or power-laws) require post-processing filtering for translating results in manufacturable data, %; unfortunately, these techniques
 which may  result case- and user-dependent and might lead to undesirable effects, like an artificial increase of the final volume fraction \citep{Zegard2016}.

Although promising results have been obtained on phase-field  topology optimization, there are still some key aspects which merit further  investigation.  Firstly, the adopted values of phase-field-related parameters affect both the obtained topology distribution and the convergence behavior of the numerical solution; hence, selecting proper values of these parameters for the specific problem at hand requires a time-consuming preliminary set of investigations. 

%REMOVED IF THE PARAMETRIC STUDY IS REMOVED: Clearly, a better insight on the influence of phase-field parameters on the final solution would support this stage, moving towards systematic optimality criteria and consistent \emph{a priori} preliminary evaluations.

Furthermore, traditional schemes for topology optimization problems generally fix the amount of material to be distributed in the design region, amount which acts as a global constraint in the optimization process.
This gold-standard strategy is easy to implement, since it adds only a single global equation in the computation; 
however, when the evolution of the topology is itself the solution of a finite element computation, the presence of such a global constraint requires to introduce global variables, altering the intrinsic element-based subdivision of the spatial discretization and introducing a coupling term between all elements. 
As an obvious consequence, the usual banded structure of the tangent stiffness matrix is altered and such an increase of the band-width of the tangent stiffness matrix might have consequences on the convergence rate of topology optimization procedures.  Currently, the literature is missing a formulation which involves only local constraints for enforcing/minimizing the amount of material to be distributed.

An additional issue of high interest is related to classically adopted numerical solution schemes. These are generally based on staggered algorithms: the evolution of the phase-field variable is obtained without the respect of the equilibrium, which is restored at the successive incremental step wherein the phase-field variable is instead frozen. In the framework of topology optimization, such an approach is referred to as Nested Analysis and Design (NAND). Allowing for the use of existing FE packages, NAND  is easy and convenient to implement and it is often preferred also as a consequence of the presence of the global constraint. Nevertheless, NAND can be computationally demanding, especially in the presence of non-linearities and for large-scale problems. On the other hand, a Simultaneous Analysis and Design (SAND) approach treats the optimization design (e.g., phase field) and equilibrium (e.g., displacement) unknowns simultaneously. In other words, the combined optimization and equilibrium problem is solved in a monolithic way by considering the evolution of topology and equilibrium conditions in a single optimization routine. Despite the fact that SAND may be competitive with respect to conventional NAND schemes, to the best of authors' knowledge, SAND (monolithic) approaches are still largely to be explored, as compared with the more classical NAND (staggered) solution strategies.

Finally, the phase{\blue -}field variable should be physically bounded in the range $[0,1]$.  To enforce this constraint, the phase-field variable obtained from the  solution of the topology evolution equation (which possibly lies outside the admissible bound) is generally projected on the admissible set at each step of the iterative solution procedure. This strategy obviously alters stationary. A variational formulation allowing for the direct enforcement of the physical bounds, without any \emph{ad hoc} algorithmic solution, would be more consistent.

Accordingly, in the described context of literature and open problems, the goals of the present paper are the following ones.
\begin{itemize}

\item To propose a mixed variational principle combining a phase-field functional for topology optimization with a (three-field) Hu-Washizu functional, then including directly in the formulation equilibrium, constitutive, and compatibility equations. The variational formulations includes also a penalization term  to directly enforce that the phase-field variable lies within the admissible physical bound $[0,1]$;

\item To specialize the mixed variational functional to a classical formulation for the topology optimization problem, referred to as {\em formulation with volume constraint}. Here, the amount of material to be distributed within the design domain is imposed \emph{a priori}, acting as a global constraint and operating on the global stiffness matrix of the corresponding finite element formulation;

\item To propose a novel formulation for the topology optimization problem, leading at the same time to a  {\em volume minimization}. The amount of material is not imposed \emph{a priori} by means of a global constraint, but rather related to a {\it cost} of the material, which is minimized together with the structural compliance. Such an approach avoids the introduction of a global constraint, respecting the convenient local nature of the finite element discretization;

\item To implement the proposed variational formulations in the context of mixed finite elements. Thanks to the combined phase-field and Hu-Washizu rationale,  a SAND monolithic algorithm is proposed, allowing   to compute the evolution law of the phase-field variable under the respect of equilibrium;

%\item To conduct, for both classical and novel formulations, a careful numerical investigation on the role of the involved parameters, highlighting the effect of their choices on the obtained solution, e.g., in terms of final compliance and computational cost;

\item To trace general guidelines for fixing the values of phase-field-related simulation parameters. In addition, a theoretical estimate for setting the {\it cost} of the material (which plays the role of a mass penalty parameter) in the volume minimization formulation as function of a target desired solution is proposed and verified;

\item To perform a comparative analysis in terms of both computational efficiency and obtained final designs for monolithic and staggered solution strategies as well as for the two investigated formulations.

\item To show the effectiveness of the proposed volume minimization formulation and its finite element mixed implementation both in two-dimensional and three-dimensional applications.

\end{itemize}

%% file: opt_cont_r1_v12.tex
\section{Theory}

\subsection{Introductory settings}

We indicate the \underline{design region} with $\Omega \subset \real^n$ (with $n=2$ or $n=3$), characterized by a volume $V$, defined as:
\begin{equation}
V= \int_{\Omega} \: d \Omega \, .
\end{equation}

The material distribution within $\Omega$ is described by a \underline{scalar phase variable} $\phi$, ideally a binary field (i.e., either 0 or $1$), with $\phi= 0$ corresponding to the absence of material (i.e., void or no material), and $\phi =1$ corresponding to the presence of material.  After topology optimization, the  material distributed inside the design region $\Omega$ occupies a  (dimensionless) volume fraction $v$ of the total volume $V$, defined as:
\begin{equation}
v= \frac{1}{V} \int_{\Omega} \phi \: d \Omega\, .
\label{eq:material}
\end{equation}

However, since a binary (i.e., either 0 or $1$) phase parameter
$\phi$ would call for sharp interfaces (between the region with material, i.e., with $\phi=1$, and the region with no material, i.e., with $\phi =0$) and since problems with  infinitely sharp interfaces are in general difficult to numerically treat and solve,  as classically done, we consider the phase parameter $\phi$ as a continuous real value field, with values in the interval $[ 0 , 1]$. Ideally, the optimal phase variable $\phi$ should be  close to an ideal binary field.

\subsection{Elastic problem}

Limiting the discussion to the case of linear elastic materials,
given a material distribution $\phi: \Omega \to [0,1]$, the \underline{elastic problem} (equilibrium, constitutive, and compatibility field equations, in combination with boundary conditions) can be expressed as:
\begin{equation}
\left\{
\begin{aligned}
& \diver \ffs + \phi \bbb = \bbzero 	& \htext{5mm}{in} & \Omega \\
& \ffs =  \fourC ( \phi) \ffe 					& \htext{5mm}{in} & \Omega \\
& \ffe = \gradsu 						& \htext{5mm}{in} & \Omega \\
& \bbu= \bbzero 						& \htext{5mm}{on} & \Gamma_D \\
& \bbt = \ffs \bbn		 				& \htext{5mm}{on} & \Gamma_N\, ,
\end{aligned}
\right.
\end{equation}
with 
$\ffs$ the stress tensor field, 
$\bbb$ the applied body-force density vector field  per unit
volume,  
$\fourC (\phi) $ a positive-definite material elastic tensor function of the phase variable $\phi$,
$\ffe$ the strain tensor field, 
$\bbu$ the displacement  vector field, 
$\diver$ and $\grads$ the divergence and the symmetric gradient operators, 
$\Gamma_D$ the region of $\partial \Omega$ where we apply Dirichlet
boundary conditions  (assumed to be  homogeneous),
$\Gamma_N$ the region of $\partial \Omega$ where we apply Neumann
boundary conditions with $\bbt$ the surface load vector and 
$\bbn$ the outward-pointing versor normal to $\Gamma_N$. 

The
elastic problem can be classically condensed in just one field equation in terms of displacements (referred to as Navier's equation),  as follows:
\begin{equation}
\left\{
\begin{aligned}
& \diver \left[ \fourC (\phi)  \gradsu  \right] + \phi \bbb = \bbzero & \htext{5mm}{in} & \Omega \\
& \bbu= \bbzero & \htext{5mm}{on} & \Gamma_D \\
& \bbt = \fourC (\phi) \gradsu  \,\bbn & \htext{5mm}{on} & \Gamma_N\, .
\end{aligned}
\right.
\label{eq:equilibrium}
\end{equation}

Assuming  voids to be  modeled as a very soft material, 
in the present work we adopt the following expression for $\fourC(\phi)$:
\begin{equation}
\fourC(\phi)=\left[\delta + (1-\delta) \frac{\text{Exp}(p \phi^p)}{\text{Exp}(p)}\right] \fourC_{A} \,  , 
\label{eq:fourC}
\end{equation}
where $\fourC_A$ is a constant positive-definite material elastic tensor (assumed to describe an isotropic behaviour with Young's modulus $E_A$ and Poisson's ratio $\nu_A$) corresponding to the solid dense material, % of the dense material, 
$p$ can be any positive value, and $0 < \delta \ll 1$ governs the low (but non-null) stiffness of the voids.  Equation \eqref{eq:fourC} is such that $\fourC(1)=\fourC_A$ and  $\fourC(0)\rightarrow \delta \fourC_A$ for $p \rightarrow +\infty$; as an example, $p =10$ already corresponds to $\|\fourC(0)- \delta \fourC_A\|<10^{-4}$. 

\begin{figure}[tbh]
\centering
\includegraphics[width=0.85\textwidth]{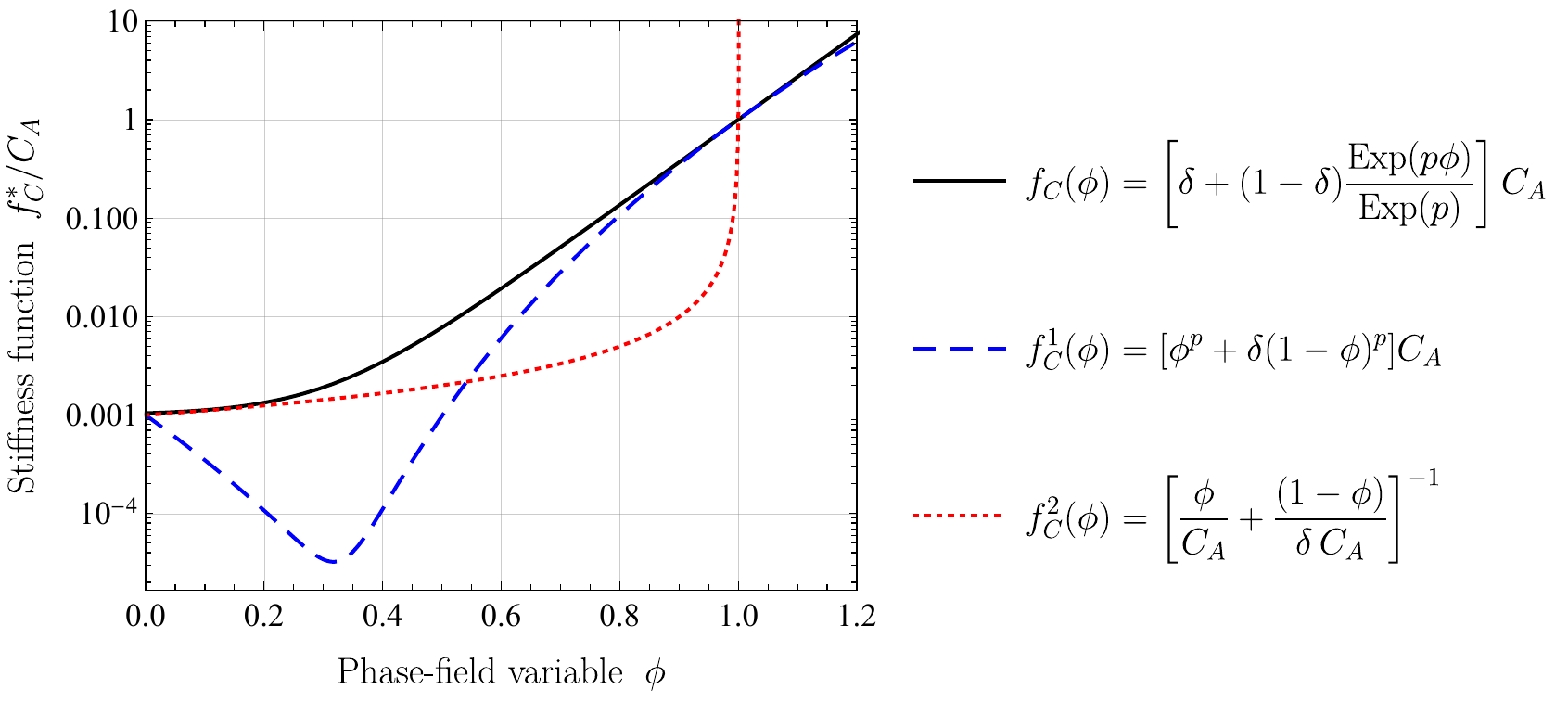}
        \caption{
        Stiffness function laws: $f_C^1(\phi)$, scalar counterpart of Eq. \eqref{eq:fourC_reg1_1}; $f_C^2(\phi)$, scalar counterpart of Eq. \eqref{eq:fourC_reg1_2}; and $f_C(\phi)$, scalar counterpart of the proposed expression \eqref{eq:fourC}. Parameters: $C_A=10$, $\delta=10^{-3}$ and $p=10$. }
        \label{fig:stiffness_function}
\end{figure}

It is worth highlighting that the particular choice of $\fourC(\phi)$ from Equation
\eqref{eq:fourC} is different from state-of-the-art expressions, which generally might read as one of the following:
\begin{subequations} \label{eq:fourC_reg1}
\begin{align}
& \fourC^1(\phi) = \left[ \phi^p + \delta (1 - \phi)^p \right] \fourC_A \,  , \label{eq:fourC_reg1_1} \\
& \fourC^2(\phi) = \left[ \phi \fourC_A^{-1} + (1 - \phi) \left(\delta \fourC_A\right)^{-1}  \right]^{-1}\, . \label{eq:fourC_reg1_2}
\end{align}
\end{subequations}
As shown in Fig. \ref{fig:stiffness_function} for their scalar counterparts, these lasts expressions either bring to non-monotonic expressions (i.e., $\fourC^1(\phi)$) or fast diverge as soon as $\phi$ exceeds the admissible range $[0,1]$ (i.e., $\fourC^2(\phi)$). 
On the one hand, a monotonically increasing behaviour would be
beneficial in terms of stability of the ensuing numerical scheme and
it would be more consistent from the physical point of view. On the
other hand, a smooth behaviour outside the admissible range would be
beneficial because the property $\phi\in[0,1]$ is sometimes violated
in numerical iterative solution schemes, often requiring an \emph{ad
  hoc} truncation of the stiffness matrix or of the phase-field
variable, and hence introducing  non-smoothness effects.  

The  expression proposed in Eq. \eqref{eq:fourC} and adopted in this
work allows to overcome drawbacks inherited by definitions in
Eq. \eqref{eq:fourC_reg1} (see Fig. \ref{fig:stiffness_function}). In
fact, $\fourC(\phi)$  turns out to be  monotonically
increasing,  everywhere-defined,  and uniformly
positive-definite for all $\phi$, being at the same time smooth with
all its derivatives.  In addition, our choice for  $\fourC(\phi)$
entails that the elastic energy density $(\phi,\ffe)\mapsto
\fourC(\phi)\ffe:\ffe/2$ is convex, contributing
to the stability of the approximation. 

\subsection{Topology optimization}

The \underline{goal of topology optimization} is  then to design an
optimal structure, i.e., to identify the optimal material distribution
$\phi$, such that structural compliance is minimized when  solving
 the elastic problem. We introduce the \underline{structural compliance} $\cC$, with the unit of measure of a work, as:
\begin{equation}
\cC (\phi, \bbu) = \int_{\Omega} \phi \bbb \cdot \bbu \: d \Omega 
+ \int_{\Gamma_N} \bbt \cdot \bbu \: d \Gamma\, ,
\end{equation}
where $\bbb$ and $\bbt$ are the actual body and surface loads (assumed to be constant), while $\bbu$ is the displacement, satisfying the elastic problem \eqref{eq:equilibrium}.

Accordingly, the \underline{topology optimization problem} can be written as:
\begin{equation}
{\text{min}_\phi} \{\cC (\phi, \bbu) \ : \ \bbu \ \text{solves} \ \
  \eqref{eq:equilibrium} \  \text{given} \ \phi\} \, .
\end{equation}
It  should be remarked  that in \eqref{eq:equilibrium} the
boundary conditions are defined on a fixed subset of $\partial
\Omega$.  Once a solution is numerically obtained, one has to
check whether $\phi\not=0$ at least on some portion of $\Gamma_D$ and
$\Gamma_N$. This generally follows from minimality. In particular, it is fulfilled by
all our computations. 

\subsection{Phase-field topology optimization}

As mentioned, the phase parameter $\phi$ is allowed to take values in the interval $[ 0 , 1]$ and, accordingly,
the region where $0 < \phi < 1$ corresponds to a smooth material interface, representing a transition from void to solid material. 
To introduce a measure of the extension or the thickness of such a transition interface and to favor fields taking values as close as possible to $0$ or $1$, we introduce a new functional $\cP$ to be minimized, defined as:
\begin{equation}
\cP ( \phi ) = \int_{\Omega} \left[  \frac{\gammaphi}{2}  \| \nabla \phi \|^2 + \frac{1}{\gammaphi} \psi_0 (\phi) \right] \: d \Omega \, ,
\end{equation}
where $\gammaphi>0$ is the \underline{thickness penalization parameter} (with the unit of measure of a length) and
$\psi_0$ is the double-well potential function, defined as:
\begin{equation}
\psi_0 (\phi) = \left[ \phi \left( \phi -1 \right) \right]^2\, .
\end{equation}
For small values of $\gammaphi$, the minimization of  $\cP (\phi )$ penalizes oscillatory material distributions and favours continuous fields $\phi$ taking values close to $0$ or $1$. For $\gammaphi \to 0$, the quantity $\cP ( \phi )$ converges in the sense of $\Gamma$-limits to the perimeter functional (cf.~\cite{modica_87}), namely a measure of the \underline{perimeter} of the interfaces between regions with material ($\phi=1$) and regions with no material ($\phi=0$). 

In the following we augment the compliance functional $ \cC$ by a multiple of $\cP$, introducing the  functional $\cJ$, defined as:
\begin{equation}
\cJ (\phi, \bbu) = \cC (\phi, \bbu) + \kphi \cP (\phi ),
\label{eq:functional_J_new}
\end{equation}
where the \underline{perimeter stiffness parameter} $\kphi>0$ (force divided by length) measures the mutual relevance of structural compliance and perimeter in the minimization. 

Accordingly, the \underline{phase-field topology optimization problem} can be written as:
  \begin{equation}
    \label{eq:to2}
     {\text{min}_\phi} \{\cJ (\phi, \bbu) \ : \ \bbu \ \text{solves} \ \
  \eqref{eq:equilibrium} \  \text{given} \ \phi\} \, , 
  \end{equation}

\subsection{Phase-field constraint}

We start recalling that ideally the phase-field variable $\phi$ should
be as much as possible  a binary field but also that any value outside
the interval, i.e., $\phi \notin [ 0 , 1 ]$, is physically 
meaningless  and it can be the source of numerical instabilities.
The minimisation of $\cJ$  in Eq.  \eqref{eq:to2}  already  favors  $\phi$ to be close to 0 or 1 through the presence of $\psi_0$, but $\phi\in[0,1]$ just in the limit as $\gammaphi\to 0$;
for all $\gammaphi>0$ (which is the case of any numerically-obtained solution) the constraint $\phi\in[0,1]$ does not necessarily hold through the minimisation of $\cJ$.

With these concepts in mind, it is clear why in most iterative schemes available in the topology optimization literature the variable $\phi$ is projected onto the admissible range $[0,1]$ at each step, 
however violating the stationary conditions of the corresponding adopted functional. 

As an alternative, we propose to further augment the functional $\cJ$ by a contribution $\mathcal{B}$, defined as:
\begin{equation}
\mathcal{B}(\phi) = \int_{\Omega} \, b(\phi) d\Omega \, , \quad \text{ with } \; b(\phi) = \begin{cases} \dfrac{(\phi - 1)^2}{2} & \phi > 1 \\
0 & 0\leq \phi \leq 1\\
\dfrac{\phi^2}{2} & \phi < 0
\end{cases}\, ,
\end{equation}
and modulated by a \underline{bounding stiffness parameter} $\kappa_b>0$ (force per unit area).

\subsection{Formulation with volume constraint}

As classically found in the literature, we may now constrain the topology optimization problem so that the portion $v$ of material to be distributed inside the design domain $\Omega$ is equal to an {\it a priori} given (dimensionless) material fraction $\overline{v}\in [ 0,1]$, defined as:
\begin{equation}
v = \frac{1}{V} \int_{\Omega} \phi \: d \Omega = \overline{v}.
\end{equation}
The corresponding topology optimization problem reads as:
\begin{equation}
    \label{eq:to3}
     {\text{min}_\phi} \{\cJ (\phi, \bbu) +\kappa_b{\mathcal B}(\phi)\ : \ \bbu \ \text{solves} \ \
  \eqref{eq:equilibrium} \  \text{given} \ \phi \  \text{and} \
  v=\overline{v}\}.
\end{equation}

Before moving on, we wish to remark that problem \eqref{eq:to3} can be
equivalently rewritten as 
\begin{align}
  \nonumber
  &   {\text{min}_\phi} \Bigg\{\int_\Omega \fourC(\phi)
       \gradsu:\gradsu+\kappa_\phi {\mathcal P}(\phi) +\kappa_b{\mathcal
  B}(\phi)+\frac{\kappa_v}{2}\int_\Omega \phi^2\:d \Omega\ : \\
  &\qquad \qquad \qquad \qquad 
  \bbu \ \text{solves} \ \
  \eqref{eq:equilibrium} \  \text{given} \ \phi \ \text{and} \
  v=\overline{v}\Bigg\}\, .  \label{eq:to32}
\end{align}
Our specific choice \eqref{eq:fourC} for $\fourC$ makes the functional under
minimization to be convex, up  to the lower-order $\psi_0$-term in $ {\mathcal P}(\phi)$.

 We  propose to solve  the above minimization problems   by investigating the stationarity of the
Lagrangian $\cL^{vc}$, defined as:
\begin{align}
\cL^{vc} (\phi, \lambda, \bbu, \ffs, \ffe ) = & \cJ (\phi, \bbu ) +
                                              \kappa_b\mathcal{B}(\phi) -  \cE^{el}(\ffe,\phi) +
                                              \nonumber\\
& +\int_{\Omega} \ffs :
                                              (\ffe - \nabla^s \bbu)
                                              \: d \Omega   +
\lambda \left[ \int_{\Omega}  (\phi - \overline{v} )\: d \Omega \right],
\label{eq:functional_volume_constraint}
\end{align}
where $\cE^{el}$ is the elastic energy, defined as:
\begin{equation}
 \cE^{el}(\ffe,\phi) = \frac{1}{2} \int_{\Omega} \fourC(\phi) \ffe : \ffe\, d\Omega.
 \end{equation}
The superscript $vc$ in $\cL^{vc}$ recalls that we are here considering a \underline{v}olume-\underline{c}onstrained formulation.

\subsection{Formulation with volume minimization}

Instead of prescribing {\em a priori}  a  fraction of material $\overline{v}$ to be distributed in the design domain $\Omega$, we may aim at exploring what is the minimum amount of material $v$ we could distribute. In terms of minimization, a corresponding topology optimization problem would read
\begin{equation}
    \label{eq:to4}
     {\text{min}_\phi} \left\{\cJ (\phi, \bbu) +\kappa_b{\mathcal
       B}(\phi)+\frac{\kappa_v}{2}\int_\Omega \phi^2\:d \Omega\ : \ \bbu \ \text{solves} \ \
  \eqref{eq:equilibrium} \  \text{given} \ \phi \right\}\, ,
  \end{equation}
where $\kappa_v>0$ is a \underline{volume penalty parameter} (force per unit area), representing for example a measure of the \lq\lq{}cost\rq\rq{} of the material per unit volume. 

It is interesting to emphasize that 
in \eqref {eq:to4} we are not minimizing the amount of  material $v$
to be distributed, but rather $\int_\Omega \phi^2\:d \Omega$; this
choice makes the problem more stable and it can be proved to be
equivalent to minimize $v$, as long as $\phi$ exclusively takes values
 $0$ or $1$.
As for the volume-constrained case,  problem \eqref{eq:to4} can be
equivalently rewritten as 
\begin{align}
   \nonumber
  &   {\text{min}_\phi} \Bigg\{\int_\Omega \fourC(\phi)
       \gradsu:\gradsu+\kappa_\phi {\mathcal P}(\phi) +\kappa_b{\mathcal
  B}(\phi)+\frac{\kappa_v}{2}\int_\Omega \phi^2\:d \Omega\ : \\
  & \qquad \qquad \qquad 
  \bbu \ \text{solves} \ \
  \eqref{eq:equilibrium} \  \text{given} \ \phi \Bigg\}\, .  \label{eq:to42}
\end{align}
Again, the choice \eqref{eq:fourC} for $\fourC$ makes the functional under
minimization to be convex, up  to the lower-order $\psi_0$-term in $ {\mathcal P}(\phi)$.

We propose to solve the topology optimization problem  \eqref {eq:to4}  by looking at the stationarity of the Lagrangian $\cL^{vm}$, defined as:
\begin{equation}
\cL^{vm} (\phi, \bbu, \ffs, \ffe  ) = \cJ (\phi, \bbu ) + \kappa_b\mathcal{B}(\phi) - \cE^{el}(\ffe,\phi) +\int_{\Omega} \ffs : (\ffe - \nabla^s \bbu) \: d \Omega +
\frac{\kappa_v}{2}
 \int_{\Omega}  \phi^2 \: d \Omega.
\label{eq:functional_minimum_volume}
\end{equation}
The superscript $vm$ in $\cL^{vm}$  recalls that we are looking at a  \underline{v}olume-\underline{m}inimization formulation.
%
%\subsection{Relation between the two formulations}
%
%\comm{Michele: for me useless. To be removed}
%
%Before moving on, we wish to comment on the fact that  the volume-constrained problem in Eq. \eqref{eq:to3} and the volume-minimization problem in Eq. \eqref{eq:to4} may be related via the following nested minimization procedure
%%
%\begin{align*}
% &  {\text{min}_\phi} \left\{\cJ (\phi, \bbu) +{\mathcal
%       B}(\phi)+\frac{\kappa_v}{2}\int_\Omega \phi^2\:d \Omega\ : \ \bbu \ \text{solves} \ \
%  \eqref{eq:equilibrium} \  \text{given} \ \phi \right\} = \\
%&\quad=
%  \min_{\overline{v}} \left(\text{min}_\phi \left\{\cJ (\phi, \bbu)
%  +{\mathcal B}(\phi) +\frac{\kappa_v}{2}\int_\Omega \phi^2\:d \Omega\ : \ \bbu \ \text{solves} \ \
%  \eqref{eq:equilibrium} \  \text{given} \ \phi \  \text{and} \
%  v=\overline{v}\right\}\right)	, .
%\end{align*}
%
%To be precise, the nested minimization in the above right-hand side does not exactly correspond to Eq. \eqref{eq:to3}, since the functional is augmented by the quadratic term $(\kappa_v/2)\int_\Omega \phi^2\:d \Omega$. This discrepancy is however difficult to detect numerically. Indeed, if $\phi$ takes only values $0$ and $1$, as it would be the case for $\gammaphi\to 0$, the above extra term reads $(\kappa_v/2) V \overline{v}$, i.e., it is independent of $\phi$. In this case, the nested minimization is exactly Eq. \eqref{eq:to3}.
%

\section{Stationarity conditions and algorithmic consideration}

\subsection{Formulation with volume constraint}
\label{Formulation with volume constraint}

To compute a stationary point of the Lagrangian $\cL^{vc}$ \eqref{eq:functional_volume_constraint}, we set its variations to zero, i.e.:
\begin{equation}
\left\{
\begin{aligned}
D_{\delta \phi} \cL^{vc}& =  
\int_{\Omega} \delta \phi \bbb \cdot \bbu \: d \Omega  +
\kphi \int_{\Omega} \left[ \gammaphi  \nabla \phi \cdot \nabla \delta
  \phi  + \frac{1}{\gammaphi} \frac{\partial \psi_0}{\partial \phi}
  \delta \phi \right] \: d \Omega  + \kappa_b \int_{\Omega} \frac{\partial b}{\partial \phi} \delta \phi \: d \Omega
\\
& - \mathcal{E}^{el}_{,\phi} (\phi, \ffe) % \frac{1}{2} \int_{\Omega} \left[ \fourC^\prime (\phi)  \ffe : \ffe \right] \delta \phi \: d \Omega  
+
\lambda  \int_{\Omega}  \delta \phi \: d \Omega = 0 
\\
D_{\delta \lambda} \cL^{vc}& = 	\delta \lambda \left[  \int_{\Omega}  (\phi - \overline{v} ) \: d \Omega \right] = 0
\\
D_{\delta \bbu} \cL^{vc} & = \int_{\Omega} \delta \bbu \cdot \phi \bbb \: d \Omega 
	+ \int_{\Gamma_N} \delta \bbu \cdot \bbt \: d \Gamma
	- \int_{\Omega}  \nabla^s \delta \bbu : \ffs \: d \Omega = 0 \\
D_{\delta \ffs}  \cL^{vc} & = \int_{\Omega} \delta \ffs : \left( \ffe - \nabla^s \bbu \right) \: d \Omega = 0 \\
D_{\delta \ffe} \cL^{vc}& =  \int_{\Omega} \delta \ffe : \left(  \ffs - \fourC(\phi) \ffe  \right)  d  \Omega = 0 \, ,
\end{aligned}
\right.
\label{eq:Dfunctional_volume_constraint}
\end{equation}
where
\begin{equation} \label{eq:E_el_phi}
\mathcal{E}^{el}_{,\phi} (\phi, \ffe)= \int_{\Omega} e^{el}_{,\phi}(\phi, \ffe)\, \delta \phi  \: d \Omega \quad \text{with} \quad 
e^{el}_{,\phi}(\phi, \ffe)= \frac{1}{2} \fourC^{\prime}(\phi)  \ffe : \ffe \; \text{ and }  \; \fourC^{\prime}(\phi) = \frac{\partial \fourC}{\partial \phi} \, .
\end{equation}
We recall that the last three relations in Eq. \eqref{eq:Dfunctional_volume_constraint} are, respectively, the equilibrium, the compatibility, and the constitutive equations, which are then automatically enforced through the stationarity of $\cL^{vc}$.

Equations \eqref{eq:Dfunctional_volume_constraint} can be rewritten in a residual form as:
\begin{equation}
\bbR^{vc} (\phi, \lambda, \bbu,\ffs,\ffe) = \bbzero\, .
\label{eq:Rvc}
\end{equation}
To solve such a residual system using an iterative scheme,  a good starting point could be to set $\phi (\bbX) = c \le 1 \; \forall \bbX \in \Omega$, i.e.,  to consider a design domain homogeneously covered by a partially dense material. Natural choices  are $c=0.5$ or $c=\overline{v}$.

In addition, one can follow a classical \underline{Allen-Cahn-like approach} by introducing  a gradient-flow dynamics. To do so, the phase-field  variable $\phi$ is  assumed to be dependent on a pseudo-time variable $t$, i.e., $\phi = \phi ( \bbX , t)$, and the new Allen-Cahn volume-constraint Lagrangian $\cL^{vc}_{AC}$ is introduced:
\begin{equation}
\cL^{vc}_{AC} = 
 \frac{\tilde{\tau}_{\phi}}{2}
\int_{\Omega}    \dot{\phi}^2 \: d \Omega
+ \cL^{vc} \, ,
\label{eq:functional_volume_constraint_AC}
\end{equation}
where the parameter $\tilde{\tau}_{\phi}>0$ (unit of measure of force times squared time per unit area) represents viscosity with respect to the pseudo-time $t$. Exploring the stationary condition of $\cL^{vc}_{AC}$ with respect to $\phi$ (assuming $\delta \phi$ to be independent from $t$) returns  the condition:
\begin{equation}
D_{\delta \phi} \cL^{vc}_{AC} = 
\tilde{\tau}_{\phi}\int_{\Omega}  \dot{\phi} \,\delta \dot{\phi} \: d \Omega
+ D_{\delta \phi} \cL^{vc} = 0 \, .
\end{equation}
At the time-discrete level, denoting with a subscript $n$ quantities evaluated at the previous time $t_n$ and with no subscript quantities evaluated at the current time $t_{n+1}$, we obtain:
\begin{equation} \label{eq:discreteL_vc_AC}
D_{\delta \phi} \cL^{vc}_{AC} = 
\tau_{\phi}\int_{\Omega}  \frac{\phi - \phi_n}{\Delta t} \delta \phi \: d \Omega
+ D_{\delta \phi} \cL^{vc} = 0\, ,
\end{equation}
with $\Delta t = t_{n+1}-t_n$ and the discrete viscosity parameter $\tau_{\phi} = \tilde{\tau}_{\phi} / \Delta t$ (unit of measure of force times time per unit area) introduced such to have the discretized temporal derivative in the evolution equation \eqref{eq:discreteL_vc_AC}. In a more explicit format,
\begin{equation}
\begin{aligned}
D_{\delta \phi} \cL^{vc}_{AC} = &  
\tau_{\phi}\int_{\Omega}  \frac{\phi - \phi_n}{\Delta t} \delta \phi \: d \Omega + \\
& + 
\int_{\Omega} \delta \phi \bbb \cdot \bbu \: d \Omega +
\kphi \int_{\Omega} \left[ \gammaphi  \nabla \phi \cdot \nabla \delta
  \phi  + \frac{1}{\gammaphi} \frac{\partial \psi_0}{\partial \phi}
  \delta \phi \right] \: d \Omega  + \kappa_b\int_{\Omega}  \frac{\partial b}{\partial \phi} \delta \phi \: d \Omega + \\
& -  \frac{1}{2} \int_{\Omega} \left[ \fourC^\prime (\phi)  \ffe : \ffe \right] \delta \phi \: d \Omega     + 
\lambda  \int_{\Omega}  \delta \phi \: d \Omega = 0.
\end{aligned}
\label{eq:vcAC}
\end{equation}
We may note that the Allen-Cahn term contributes to the monotonicity of the map: 
$$\phi \mapsto D_{\delta \phi}\cL^{vc}_{AC}(\phi),$$ 
delivering enhanced stability to the numerical scheme,  although
 the map $\phi \mapsto   D_{\delta \phi}\cL^{vc}_{AC}(\phi)$. 
Indeed, such map turns out to be monotone, apart from the polynomial
  $\psi_0$-term and the bilinear term $\phi\bbb \cdot {\boldsymbol u}$
  in the compliance.  By choosing $\tau_\phi/\Delta t$ large with
  respect to $1/\gammaphi$, the Allen-Cahn term dominates the 
  $\psi_0$-term. In particular,   
the function $\phi \mapsto (\tau_\phi/2\Delta t)\phi^2 + (\kappa_\phi/\gamma_\phi) (\phi(1-\phi))^2$ is convex 
whenever
\begin{equation}
\frac{\tau_\phi}{\Delta t}\geq 2 \frac{\kappa_\phi}{\gammaphi} \, . \label{eq:tau}
\end{equation}
In the literature it is often chosen $\tau_\phi=\gammaphi$ (modulo fixing 
dimensions). This choice is of course compatible with Eq. \eqref{eq:tau}, as long as $\Delta t$ is taken to be small enough. Nevertheless, motivated by  Eq. \eqref{eq:tau}, the following relationship between $\tau_\phi$ and other parameters in the functional will be introduced:
\begin{equation}
\tau_\phi = \frac{\kappa_\phi T_{\phi}}{\gammaphi} \, , \label{eq:tau_def}
\end{equation}
where $T_{\phi}$ is a characteristic time constant. %Note that Eq. \eqref{eq:tau} is sufficient but not necessary for the convexity of $\cL^{vc}_{AC}$, since the gradient term in $\phi$ also contributes. 

Substituting Equation \eqref{eq:Dfunctional_volume_constraint}$_1$ with \eqref{eq:vcAC} we obtain the new residual equation corresponding to the formulation with volume constraint in an Allen-Cahn approach
\begin{equation}
\bbR^{vc}_{AC} (\phi, \bbu,\ffs,\ffe, \lambda) = \bbzero \, .
\label{eq:RvcAC}
\end{equation}

\subsection{Formulation with volume minimization} \label{Formulation with volume minimization}

To compute a stationarity point of the Lagrangian $\cL^{vm}$  \eqref{eq:functional_minimum_volume}, 
we set its variations to zero, i.e.:
\begin{equation}
\left\{
\begin{aligned}
D_{\delta \phi} \cL^{vm}& =  
\int_{\Omega} \delta \phi \bbb \cdot \bbu \: d \Omega +
\kphi \int_{\Omega} \left[ \gammaphi  \nabla \phi \cdot \nabla \delta \phi  + \frac{1}{\gammaphi} \frac{\partial \psi_0}{\partial \phi} \delta \phi \right] \: d \Omega  + \kappa_b\int_{\Omega} \frac{\partial b}{\partial \phi} \delta \phi \: d \Omega
\\
& - \mathcal{E}^{el}_{,\phi} (\phi, \ffe) % \frac{1}{2} \int_{\Omega} \left[ \fourC^\prime (\phi)  \ffe : \ffe \right] \delta \phi \: d \Omega  
+ 
\kappa_v  \int_{\Omega}  \phi \delta \phi \: d \Omega = 0 
\\
D_{\delta \bbu} \cL^{vm} & = \int_{\Omega} \delta \bbu \cdot \phi \bbb \: d \Omega 
	+ \int_{\Gamma_N} \delta \bbu \cdot \bbt \: d \Gamma
	- \int_{\Omega}  \nabla^s \delta \bbu : \ffs \: d \Omega = 0 \\
D_{\delta \ffs}  \cL^{vm} & = \int_{\Omega} \delta \ffs : \left( \ffe - \nabla^s \bbu \right) \: d \Omega = 0 \\
D_{\delta \ffe} \cL^{vm}& =  \int_{\Omega} \delta \ffe : \left(  \ffs - \fourC(\phi) \ffe  \right)  d  \Omega = 0 \,  .
\end{aligned}
\right. \, .
\label{eq:Dfunctional_minimu_volume}
\end{equation}

We recall again that the last three relations in
Eq. \eqref{eq:Dfunctional_minimu_volume} are, respectively, the
equilibrium, the compatibility, and the constitutive equations, which
are then automatically enforced through the  stationarity of
$\cL^{vm}$.  The term $\mathcal{E}^{el}_{,\phi} (\phi, \ffe)$ in Eq. \eqref{eq:Dfunctional_minimu_volume}$_1$ is given in Eq. \eqref{eq:E_el_phi}.

Equations \eqref{eq:Dfunctional_minimu_volume} can also be rewritten in a residual form as:
\begin{equation}
\bbR^{vm} (\phi,  \bbu,\ffs,\ffe) = \bbzero.
\label{eq:Rvm}
\end{equation}
To solve such a residual  system using an iterative scheme, as
for the formulation with volume constraint, we can follow an \underline{Allen-Cahn-like approach}, i.e., introducing  a gradient-flow dynamics.  To do so, the phase-field  variable $\phi$ is again assumed to be dependent on a pseudo-time variable $t$, i.e., $\phi = \phi ( \bbX , t)$, and a new Lagrangian $\cL^{vm}_{AC}$ is introduced:
\begin{equation}
\cL^{vm}_{AC} = 
 \frac{\tilde{\tau}_\phi}{2}
\int_{\Omega}  \dot{\phi} ^2 \: d \Omega
+ \cL^{vm}.
\label{eq:functional_minimu_volume_AC}
\end{equation}
The stationarity of $\cL^{vm}_{AC}$ with respect to $\phi$ (assuming $\delta \phi$ to be independent from $t$) corresponds to 
\begin{equation}
D_{\delta \phi} \cL^{vm}_{AC} = 
{\tau_\phi}\int_{\Omega}  \frac{\phi - \phi_n}{\Delta t} \delta \phi \: d \Omega
+ D_{\delta \phi} \cL^{vm} = 0\, ,
\end{equation}
with $\tau_\phi = \tilde{\tau}_\phi /\Delta t$. In a more explicit format 
\begin{equation}
\begin{aligned}
D_{\delta \phi} \cL^{vm}_{AC} = &  
\tau_\phi \int_{\Omega}  \frac{\phi - \phi_n}{\Delta t} \delta \phi \: d \Omega + \\
& + 
\int_{\Omega} \delta \phi \bbb \cdot \bbu \: d \Omega +
\kphi \int_{\Omega} \left[ \gammaphi  \nabla \phi \cdot \nabla \delta \phi  + \frac{1}{\gammaphi} \frac{\partial \psi_0}{\partial \phi} \delta \phi \right] \: d \Omega  + \kappa_b\int_{\Omega} \frac{\partial b}{\partial \phi} \delta \phi \: d \Omega +
\\
& -  \frac{1}{2} \int_{\Omega} \left[ \fourC^\prime (\phi)  \ffe : \ffe \right] \delta \phi \: d \Omega   +
\kappa_v  \int_{\Omega}  \phi \delta \phi \: d \Omega = 0\, . 
\end{aligned}
\label{eq:vmAC}
\end{equation}

The Allen-Cahn term again contributes to the  monotonicity of
$\phi\mapsto D_{\delta \phi}  \cL^{vm}_{AC}(\phi)$.  Compared with the volume-constraint case, the situation is here more favorable, since the $\kappa_v$ term is also contributing to convexity. In particular, the function $\phi \mapsto (\tau_\phi/2\Delta t)\phi^2 + (\kappa_\phi/\gamma_\phi) (\phi(1-\phi))^2+(\kappa_v/2)\phi^2$ is convex whenever
\begin{equation}
\frac{\tau_\phi}{\Delta t}+\kappa_v\geq 2\frac{\kappa_\phi}{\gammaphi} \, .
 \label{eq:tau2}
\end{equation}
In particular, comparing  \eqref{eq:tau2} with \eqref{eq:tau}, one has
that from \eqref{eq:tau2} convexity still holds for  $\tau_\phi=0$  (no Allen-Cahn regularization),  as long as $\kappa_v$ is large enough. Also for this formulation, $\tau_{\phi}$ will be defined as function of other parameters as in Eq. \eqref{eq:tau_def}.

Substituting Equation \eqref{eq:Dfunctional_minimu_volume}$_1$ with \eqref{eq:vmAC} we obtain the new residual system corresponding to the formulation with volume constraint in an Allen-Cahn approach
\begin{equation}
\bbR^{vm}_{AC} (\phi, \bbu,\ffs,\ffe) = \bbzero\, .
\label{eq:RvmAC}
\end{equation}

\subsection{Solution algorithm, convergence criteria, and output quantities} \label{sec:algorithm}

In a time-discrete setting, at each time instant $t_{n+1}$ we need to solve the non-linear residual equation  \eqref{eq:RvcAC} for the formulation with volume constraint and  the non-linear residual equation \eqref{eq:RvmAC} for the formulation with volume minimization. Furthermore, assuming to properly solve at each time instant $t_{n+1}$ the corresponding non-linear residual problem, we need to follow the time-discrete dynamics to convergence to a stationary solution. In the following, stationarity is measured by controlling changes in time of the material distribution $\phi$ and of the displacement field $\bbu$.

Accordingly, we adopt an algorithm with a double iteration loop, i.e., an external iteration loop on time (controlling a norm on material distribution $\phi$ and  displacement field $\bbu$ changes)  and an internal iteration loop (controlling a norm on the satisfaction of the residual problem).

The \underline{stopping criterion for the external iteration loop on time} is defined on the basis of the following  error relative to the Allen-Cahn procedure, or briefly  \underline{Allen-Cahn error}:
\begin{equation}
\cE_{AC} =  \frac{\cE_{\Delta \phi} + \cE_{\Delta u}}{2}\, ,
\end{equation}
where:
\begin{equation}
 \cE_{\Delta \phi} = \frac{T_{\phi}}{V \Delta t }\int_{\Omega} \| \phi - \phi_n \| d\Omega\, , \quad
\cE_{\Delta u} = \frac{T_{\phi}}{A_{\Delta u} V \Delta t}\int_{\Omega} \| \bbu - \bbu_n \| d\Omega \, ,
\end{equation}
with $T_{\phi}$ the time constant introduced in Eq. \eqref{eq:tau_def} and $A_{\Delta u}$ a normalization constant, computed as:
\begin{equation}
A_{\Delta u} = \frac{1}{V} \int_{\Omega} \| \bar{\bbu}_{sol} \| d\Omega \, .
\end{equation}
Here, $\bar{\bbu}_{sol}$ is the reference displacement field obtained by solving the equilibrium problem with a fixed  $\phi = 1$ everywhere in the design domain.
The iterations in time are stopped when $\cE_{AC} < c_{AC}^{conv}$, with $c_{AC}^{conv}$  an \ul{Allen-Cahn convergence parameter} to be set. 

At each time instant $t_{n+1}$, the internal iteration loop consists
in the solution of the residual systems \eqref{eq:RvcAC} and
\eqref{eq:RvmAC} by means of the Newton-Raphson method which is
implicit since time discretization has been performed by means of a
backward Euler strategy. Starting from the functionals $\mathcal{L}$,
derivatives required for computing the residual vector $\bbR$ and the
tangent system matrix $\mathbb{D}$ (i.e., linearization of $\bbR$) are
performed by means of the Mathematica package AceGen
\citep{Korelc2016},  allowing  for combined symbolic-numeric programming. The solution of the finite element resulting system is obtained by means of the Mathematica package AceFEM \citep{Korelc2016}. The \underline{stopping criterion for the internal iteration loop} is defined on the basis of the norm of  the residual vector $\bbR$  corresponding to the governing equations, i.e.:
$$
\cE_R = \| \bbR_{n+1} \|\, .
$$
The iterations on the residual equations are stopped when $\cE_{R} < c_{R}^{res}$, with $c_{R}^{res}$ a 
\ul{residual convergence parameter} to be set, returning the updated solution of primary variables at time instant $t_{n+1}$. Divergence is obtained when a maximum number $N_{iter}^{max}=15$ of Newton-Raphson iterations is reached. In this case, the time-step size $\Delta t$ is decreased and the internal loop at $t_{n+1}$ is started again. The tuning of the time-step size is performed throughout the solution by means of a path-following procedure with an adaptive time stepping which returns the optimal $\Delta t$ around a reference value $\Delta t_0$, here chosen in the range $\Delta t \in (10^{-5}, 10^3)\Delta t_0$, \citep{Korelc2016}. 

An overview of the solution algorithm is reported in Algorithm \ref{algo:algo1}. It is noteworthy that, in our approach, we do not solve residual equations in a staggered form (as sometimes proposed in the literature) but we numerically  iterate in an implicit form on all residual equations, which means that our approach is fully implicit and monolithic.
In other words, the implemented solution strategy corresponds to a Simultaneous Analysis and Design (SAND) approach. For the sake of comparison, results will be compared also by adopting a Nested Analysis and Design (NAND) approach, which corresponds to a staggered implementation between the phase-field evolution problem and the mechanical equilibrium conditions. This is achieved by assuming that:
\begin{itemize}
\item The elastic energy rate term in Eqs. \eqref{eq:Dfunctional_volume_constraint}$_1$ and  \eqref{eq:Dfunctional_minimu_volume}$_1$ is computed at the previous time step, that is $\mathcal{E}^{el}_{,\phi}(\phi_n,\boldsymbol{\varepsilon}_n)$ (see Eq. \eqref{eq:E_el_phi});

\item The constitutive response considered in the mechanical equilibrium balance (i.e., Eqs. \eqref{eq:Dfunctional_volume_constraint}$_5$ and  \eqref{eq:Dfunctional_minimu_volume}$_4$) is computed on the basis of the phase-field computed at the previous time step, that is $\mathbb{C}(\phi_n)$.

\end{itemize}
For the staggered NAND approach, the phase-field variable is projected
within the admissible range $[0,1]$ at each step, in agreement with
existing literature \citep[e.g.][]{carraturo_19}.  This projection has been verified to be necessary for convergence issues, since otherwise $\fourC(\phi_n)$ takes values leading to instability/divergence of the numerical iteration scheme. On the other hand, for the monolithic SAND approach, the phase-field variable is not \emph{ad hoc} projected at each step. This choice is more sound from the theoretical point of view since the bounding constraint is already included in the variational functional, and obtained results will show that it does not lead to divergence issues.

%\fbox{
\begin{algorithm}[tbh]
 \KwData{Geometrical, material and simulation parameters}
Define and discretize design region $\Omega$\\
Define solution vector $\bbx = \{\phi, \bbu,\ffs, \ffe,\lambda\}$ for $\mathcal{L}^{vc}$ and $\bbx = \{\phi, \bbu,\ffs, \ffe\}$ for $\mathcal{L}^{vm}$

\While{\rm [ $\cE_{AC} > c_{AC}^{conv}$
   	\htext{2mm}{and}
	$t_{n+1} \le t_{final}$ ]
	}{
Initialize $\bbR_{n+1}=\bbR_{n}$ and $\mathbb{D}_{n+1}=\mathbb{D}_{n}$\\
	\While{\rm [ $\cE_{R} > c_{R}^{res}$ 
   	\htext{2mm}{and}
	$N_{iter} \le N_{iter}^{max}$ ] 
   	}{
   Update solution: $\bbx_{n+1} \rightarrow \bbx_{n} - \mathbb{D}_{n+1}^{-1}\bbR_{n+1}$ \\
   Update residual $\bbR_{n+1}$ and tangent $\mathbb{D}_{n+1}$ \\
   Compute residual norm $\cE_R = \| \bbR_{n+1} \|$
 }
 \eIf{$\cE_{R} \leq c_{R}^{res}$}{
   Compute  $\cE_{AC}$ \\
   Compute new $\Delta t$ from adaptive time stepping algorithm \\
   Increment time instant $t_{n+1} = t_n + \Delta t$ \\
   Update time step: $n \rightarrow n+1$ \\
   }{
   Decrease $\Delta t$ from adaptive time stepping algorithm \\
   Re-define $t_{n+1} = t_n + \Delta t$ \\
   Go to line 4 \\
  }
%Update solution $\bbx_{n+1} \rightarrow \bbx_n$
}
 \KwResult{Return $\bbx_{sol}=\bbx_{n}$ and compute post-processing quantities}
 \caption{Implicit and monolithic solution algorithm for the phase-field topology optimization problem with an Allen-Cahn-like strategy (SAND approach).}
 \label{algo:algo1}
\end{algorithm}
%}

Once a converged final solution is obtained (and in particular the corresponding phase-field $\phi_{sol}$, displacement $\bbu_{sol}$, and stress $\ffs_{sol}$), the solution quality is
expressed  by  computing and evaluating the following three quantities:
\begin{itemize}
\item
The volume fraction $v_{sol}$,  
$$
v_{sol} = \frac{1}{V}\int_{\Omega} \phi_{sol} d\Omega\, ;
$$
%\item
%the structural compliance $\cC_{sol}=\cC(\phi_{sol},\bbu_{sol})$;
\item
A scalar measure of the structural displacement $U_{sol}$, obtained by dividing the final value of the compliance $\cC_{sol}=\cC(\phi_{sol},\bbu_{sol})$ with the total resultant $F_{tot}$ of the applied loads:
$$
U_{sol} = \frac{\cC_{sol}}{F_{tot}} \qquad \text{with } \quad F_{tot}=\int_{\Omega} \phi_{sol} \bbb d\Omega + \int_{\Gamma_N} \bbt d\Gamma\;  ;
$$

\item
The diffused-perimeter measure $\mathcal{P}_{sol}=\cP(\phi_{sol})$.
\end{itemize}

To compare different solutions (in particular corresponding to different material distributions) from an engineering point of view, we also compute the Von Mises stress distribution $\sigma_{vm}({\bbX})$, 
the maximum Von Mises stress $\sigma_{vm}^{max}$, and the 
average Von Mises stress $\sigma^{avg}_{vm}$, defined respectively as:
$$
\sigma_{vm} ({\bbX}) = \| \bbs_{sol} ({\bbX}) \|,
$$
$$
\sigma_{vm}^{max}=\max_{\bbX \in \Omega} \{ \sigma_{vm}({\bbX}) \},
$$
\begin{equation}
\sigma_{vm}^{avg} = \frac{1}{V_{sol}}\int_{\Omega_f} \sigma_{vm}({\bbX}) d\Omega
\label{eq:sigma_vm}
\end{equation}
where $\bbs_{sol}$ is the deviatoric part of the stress tensor $\ffs_{sol}$ at solution, $\| \cdot \|$ is the standard Euclidean norm, $V_{sol}$ is the effective total volume occupied by the material at solution, i.e.
$$
V_{sol} = \int_{\Omega} \phi_{sol} d\Omega = v_{sol} V\,,
$$
and $\Omega_f= \{\bbX \;\text{ s.t. } \;\phi(\bbX)>0.5\}$ is the filled domain.

%% file: opt_discretization_r1_v12.tex
\section{Finite-element approximation and parameters settings}

To solve the problem  under investigation,  we introduce a space discretization, considering a standard decomposition of the design domain $\Omega$ into a set of non-overlapping finite elements, and for each single element we introduce the  approximations described in the following. Both two-dimensional (2D) and three-dimensional (3D) examples will be considered.

\begin{itemize}
\item
The displacement field $\bbu$ is approximated in a $C^0$-continuous element-by-element  form as:
\begin{equation}
\bbu =
\sum_{i=1}^{n_n^e} N_i(\ffxi) \bbuhatk{i}
\end{equation}
%\begin{equation*}
%\bbu = 
%\left\{
%\begin{array}{c}
%u \\ v
%\end{array}
%\right\} =
%\sum_{i=1}^4 N_i \bbuhatk{i}
%\end{equation*}
%
where $N_i(\ffxi)$ are Lagrange-type shape functions defined on the parent coordinate system $\ffxi$  and $ \bbuhatk{i}$ are the nodal degrees of freedom, with $n_n^e$ being the total number of nodes for each element. A standard isoparametric mapping from the parent coordinate system $\ffxi$ to the spatial (physical) coordinate system $\bbX$ is implemented.

\item
The phase field $\phi$ is approximated in a C$^0$-continuous element-by-element form as:
\begin{equation}
\phi = 
\sum_{i=1}^{n_n^e} N_i(\ffxi) \hat{\phi}_i \, ,
\end{equation}
%
%$$
%\phi = 
%\sum_{i=1}^4 N_i \hat{\phi}_i \, ,
%$$
%
where $\hat{\phi}_i$ denotes nodal degrees of freedom and $N_i$, $n_n^e$ are defined above.

\item
The stress field $\ffs$  is approximated in a discontinuous element-by-element form, introducing element unknowns collected in  vector $\hat{\ffs}$, \citep{Djoko2006}. In the parent coordinate system $\ffxi$, the second-order stress tensor $\ffs_{\ffxi}$ is interpolated as:
\begin{equation}
\text{vec}(\ffs_{\ffxi}) =
\bbN_{\ffs}(\ffxi) \hat{\ffs}\, ,
\end{equation}
where $\text{vec}(\cdot)$ denotes the vector Voigt representation of stresses, and $\bbN_{\ffs}$ represents the shape functions for the stresses. In the spatial (physical) coordinate system $\bbX$, the second-order stress tensor $\ffs$  is computed from the one in the parent domain, $\ffs_{\ffxi}$, via:
\begin{equation} \label{eq:sigma_transf}
\ffs(\bbX) = \bbT \ffs_{\ffxi}(\ffxi)  \bbT^T\, ,
\end{equation}
where $\bbT$ is a  second-order transformation tensor which allows to map stresses from the parent $\ffxi$ to the spatial $\bbX$ cartesian space. The transformation matrix must:
\begin{enumerate}
\item Produce stresses in spatial cartesian space which satisfy the patch test (i.e., can produce constant stresses and be stable);

\item Be independent of the orientation of the initially chosen element coordinate system and numbering of element nodes (invariance requirement).
\end{enumerate}

For the second-order transformation tensor $\bbT$,  \cite{PianSumihara1984} proposed a constant array (to preserve constant stresses) deduced from the Jacobian ${\bf J}(\ffxi)$ associated to the geometrical mapping between the physical coordinates $\bbX$ and the parent coordinates $\ffxi$, computed at the centre of the element, i.e., $\bbT ={\bf J}|_{\ffxi=\bbzero} = {\bf J}_0$. In the present paper, we employ the normalized counterpart $\tilde{\bf J}_0$ of it:
\begin{equation}
\bbT = \tilde{\bf J}_0= \frac{{\bf J}_0}{\sqrt{\text{Det} {\bf J}_0}} \quad \text{ with } \; {\bf J}_0= \left. \frac{\partial \bbX ( \ffxi )}{\partial \ffxi} \right|_{\ffxi= \bbzero}\, .
\end{equation}

\item
The strain field $\ffe$  is approximated in a discontinuous element-by-element form, introducing element unknowns collected in vector $\hat{\ffe}$, \citep{Djoko2006}. In the parent coordinate system $\ffxi$, the second-order strain tensor $\ffe_{\ffxi}$ is interpolated as:
\begin{equation}
\text{vec}(\ffe_{\ffxi}) =
\bbN_{\ffe}(\ffxi) \hat{\ffe} \, ,
\end{equation}
where $\bbN_{\ffe}$ represent the shape functions for the strains. In the spatial coordinate system $\bbX$, the second-order strain tensor $\ffe$ is computed from the one $\ffe_{\ffxi}$ through the transformation tensor ${\bf T}$ as\footnote{%
This choice satisfies with Eq. \eqref{eq:sigma_transf} the principle of energy equivalence in the parent and spatial coordinate systems:
$$
\int_{\Omega_e} \ffs : \ffe d \Omega =
\int_{\Omega_{\Box}} \ffs : \ffe J{\ffxi} d \Box =
\int_{\Omega_{\Box}}  \ffs_{\ffxi} : \ffe_{\ffxi} J_0 d \Box\, ,
$$
which immediately follows since
\begin{equation*}
\begin{split}
\ffs : \ffe =
\bbT \ffs_{\ffxi}  \bbT^T : \bbT^{-T} \ffe_{\ffxi} \bbT^{-1} =
\bbT \ffs_{\ffxi}  : \bbT^{-T} \ffe_{\ffxi} \bbT^{-1} \bbT=
\bbT \ffs_{\ffxi}  : \bbT^{-T} \ffe_{\ffxi}  =
\ffs_{\ffxi}  : \bbT^T \bbT^{-T} \ffe_{\ffxi}  =
\ffs_{\ffxi} : \ffe_{\ffxi}\, .
\end{split}
\end{equation*}
}:
\begin{equation}
\ffe(\bbX) = \frac{J_o}{J(\ffxi)} \bbT^{-T} \ffe_{\ffxi}(\ffxi)  \bbT^{-1}\, ,
\end{equation}
where $J(\ffxi)=\text{Det}({\bf J}(\ffxi))$ and $J_o=\text{Det}({\bf J}_o)$.

\item
For the case of the formulation with volume  constraint (i.e., $\mathcal{L}^{vc}$), the Lagrange multiplier 
$\lambda$ is assumed as a constant global field, hence:
\begin{equation}
\lambda = \hat{\lambda}\, .
\end{equation}

\end{itemize}

In numerical applications, four node quadrilateral elements (QUAD-4 or Q1 elements) are implemented for 2D applications, while eight node hexahedral elements (HEX-8 or H1 elements) for 3D applications. Therefore, Lagrange-type shape functions $N_i$ correspond to bilinear and trilinear polynomials, respectively for 2D and 3D applications. More information on the interpolation of the stress and strain fields, together with the specific forms of $\bbN_{\ffs}$ and $\bbN_{\ffe}$ employed in this work \citep{Weissman96,Cao2002,Djoko2006}, are given in Appendix \ref{app:interpolation}.  Static condensation of stress and strain variables is performed at element level to reduce the dimensions of the numerical problem, which then correspond to the one of a pure-displacement phase-field formulation.

\subsection{Settings of simulation parameters} \label{sec:par_settings}

The performance of phase field formulations highly depend on the values of the chosen parameters. These can be divided in physical and simulation parameters. 

On one hand, physical parameters  are  the material properties contained in the tangent stiffness matrix $\fourC$ and, for the volume constraint formulation $\cL^{vc}$, the final amount of material distributed in the design domain.

On the other hand, simulation parameters are the thickness penalization parameter $\gamma_{\phi}$ (with the unit of length), the perimeter stiffness parameter $\kappa_{\phi}$ (with the unit of force per unit length), and the bounding stiffness parameter $\kappa_b$ (with the unit of force per unit area). Moreover, for the volume minimization formulation $\cL^{vm}$, the final amount of material distributed in the design domain (i.e., $v_{sol}$) depends on another simulation parameter represented by the volume penalty parameter $\kappa_v$  (with the unit of force per unit area).

The following considerations are traced for properly settings the values of simulation parameters. As regards the volume constraint formulation $\cL^{vc}$:

\begin{itemize}
\item The thickness penalization parameter $\gamma_{\phi}$ is chosen comparable to a typical element dimension $h_e$ in the mesh discretization, that is $\gamma_{\phi} \approx h_e$;

\item The perimeter stiffness parameter $\kappa_{\phi}$ is chosen such that the term $\kappa_{\phi}/\gamma_{\phi}$ (governing the driving force which allows to obtain a binary-like solution for the phase field) is comparable with the strain-energy density rate $e^{el}_{,\phi}$ in Eq. \eqref{eq:E_el_phi}.  The latter quantity is clearly not constant during the solution. A proper  estimate can be anyway obtained as $\kappa_{\phi} \approx \gamma_{\phi} \bar{e}^{el}_{,\phi} $, with $\bar{e}^{el}_{,\phi}=e^{el}_{,\phi}(1,\bar{\boldsymbol{\varepsilon}})$ being the strain-energy density rate obtained on the \lq\lq{}full material case\rq\rq{}. In detail, $\bar{\boldsymbol{\varepsilon}}$ is the reference strain field obtained by solving the equilibrium problem with a fixed  $\phi = 1$ everywhere in the design domain;

\item The bounding stiffness parameter $\kappa_b$  should be  at least 3 order of magnitudes greater than $\kappa_{\phi}$, that is $\kappa_b>10^3 \kappa_{\phi}$.
 
\end{itemize}

As regards the volume minimization formulation $\cL^{vm}$, the same
rules-of-thumb can be introduced for $\gamma_{\phi}$ and
$\kappa_b$. Nevertheless, contrarily to the volume constraint
formulation, additional  care is needed in order to obtain an 
% considerations should be traced if one wants to have an
\emph{a priori} estimate of the final volume fraction $v_{sol}$. In
fact, although the volume minimization principle inherits advantages
from the design viewpoint (as the following results will show), a
preliminary estimate of $v_{sol}$ would be useful % at least
for obtaining a reference design solution. Therefore, the following is proposed:
\begin{itemize}
\item The volume penalty parameter $\kappa_v$ is determined on the basis of a target estimate of the final volume fraction $v_{sol}^{tar}$. To this aim, a function $v_{sol}^{tar}(\kappa_v)$ is introduced such that the following three conditions are met:
\begin{equation}\label{eq:v_sol^est_cond}
v_{sol}^{tar}|_{\kappa_v=0}=1 \, , \quad \lim_{\kappa_v\rightarrow +\infty} v_{sol}^{tar} = 0 \, , \quad \left.\frac{\partial v_{sol}^{tar}}{\partial \kappa_v}\right|_{\kappa_v=0} = -\frac{1}{\bar{e}^{el}_{,\phi}}\, .
\end{equation}
The first two conditions enforce the limit conditions that a material
with \lq\lq{}null cost\rq\rq{} (i.e., $\kappa_v=0$) will fill the
entire domain design, while for a material with an extremely
\lq\lq{}high cost\rq\rq{} (i.e., $\kappa_v \rightarrow +\infty$) the
domain design will tend to be empty. The last condition prescribes
that, starting from \lq\lq{}null material cost\rq\rq{}, the variation
of $v_{sol}^{tar}$ associated with an increase of $\kappa_v$  is
inversely proportional to $\bar{e}^{el}_{,\phi}$. Conditions
\eqref{eq:v_sol^est_cond} are met by the following definition of 
the   function $v_{sol}^{tar}=v_{sol}^{tar}(\kappa_v)$:
\begin{equation}\label{eq:v_sol^est}
v_{sol}^{tar}(\kappa_v) = \frac{\bar{e}^{el}_{,\phi}}{\kappa_v + \bar{e}^{el}_{,\phi}}\, ,
\end{equation}
which allows to fix $\kappa_v$ on the basis of the target desired value of $v_{sol}^{tar}$;

\item Once $\kappa_v$ is  fixed,  the perimeter stiffness
  parameter $\kappa_{\phi}$ is now chosen such that the term
  $\kappa_{\phi}/\gamma_{\phi} \approx \kappa_v$, or $\kappa_{\phi}=
  \kappa_v \gamma_{\phi}$. This choice follows the same  rationale
  of  the volume constraint formulation, but it allows to tune
  $\kappa_{\phi}$ as  a  function of the target volume fraction.
\end{itemize}

%% file: opt_results_r1_v12.tex
\newpage

\section{Numerical results}

The following Section presents an extensive numerical investigation resulting in a detailed convergence study and a discussion on the obtained final designs. The numerical results clearly highlight differences between the two investigated formulations as well as advantages related to the monolithic solution strategy, with numerical  simulations   addressing both two-dimensional and three-dimensional applications.

\begin{figure}[tbh]
\centering
\includegraphics[width=0.99\textwidth]{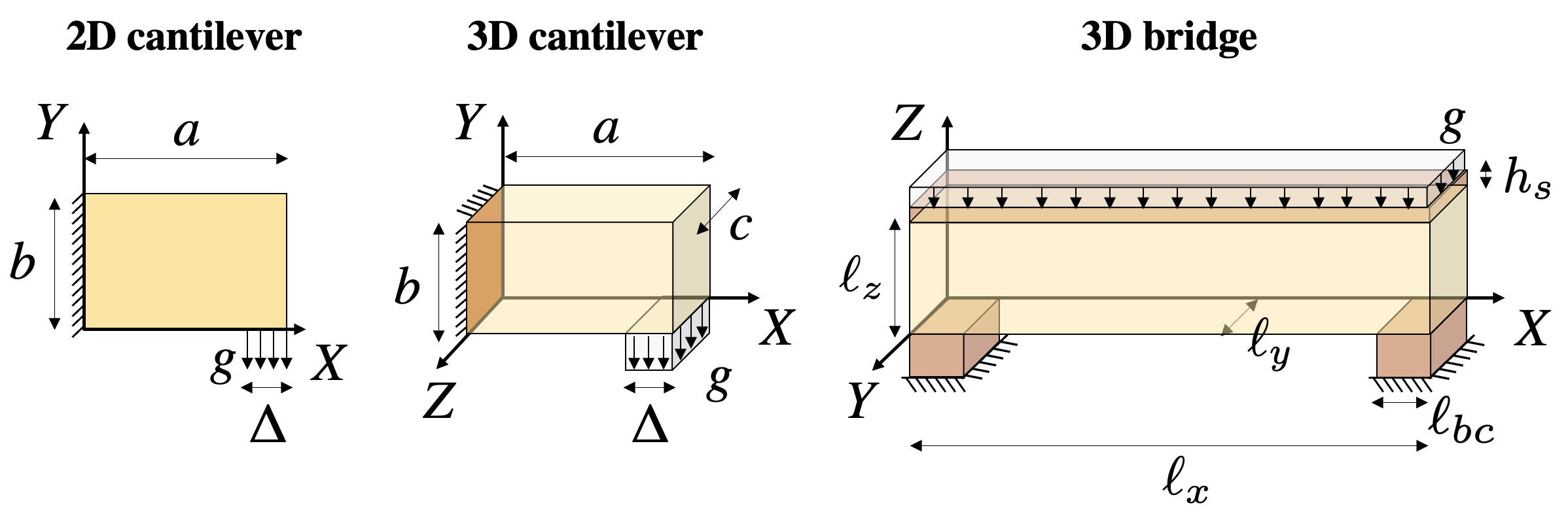}
        \caption{
        Design domain, loading and boundary conditions of numerical simulations: the 2D cantilever (Section \ref{sec:res2D}); the 3D cantilever (Section \ref{sec:res3D}); and the 3D bridge (Section~\ref{sec:res3D}).}
        \label{fig:test1_geometry}
\end{figure}

In particular, Section \ref{sec:res2D} presents the results of a 2D cantilever and  performs  a comparison between different formulations and solution approaches, while Section \ref{sec:res3D}  extends  the discussion to two 3D examples, i.e.,   a 3D cantilever  and a 3D bridge (see Figure \ref{fig:test1_geometry}).

In all the simulations, null body forces are considered, e.g., $\bbb={\bf 0}$. Moreover, the characteristic time constant is chosen as $T_{\phi}=1$ s, the discrete viscosity parameter $\tau_{\phi}$ as in Eq. \eqref{eq:tau_def}, the reference time step increment $\Delta t_0=10^{-2}$ s, and the  convergence parameter for the residual internal iteration loop fixed as $c_{R}^{res}=10^{-8}$ J.

\subsection{2D simulations: parametric analyses} \label{sec:res2D}

The 2D cantilever is solved starting from a two-dimensional design domain of rectangular shape with dimensions $a \times b$ (see Fig. \ref{fig:test1_geometry}), assuming plane strain conditions, with out-of-plane thickness $c$. 
The domain is clamped on the left vertical side and loaded by a downward vertically-oriented surface traction $g$ on  a segment of length $\Delta$ at the right bottom end of the domain. 

Table \ref{tab:par1} reports  geometrical and material parameters. If not differently specified, the initial conditions for the phase variable $\phi$ in the entire domain $\Omega$ are set as follows:  $\phi_0 = \bar{v}$ for all simulations relative to the formulation with volume constraint; $\phi_0 = 1$ for all simulations relative to the formulation with volume minimization. Colour maps of   the phase-field  variable $\phi$ are represented in the region of $\Omega$ where the phase variable $\phi \ge 0.1$. 
 Moreover, 
 %
 %compliance $\mathcal{C}$ will be presented as normalized with respect to the total applied force magnitude $F_g = g \times c \times \Delta$, while 
% 
 the simulation time is normalized with respect to the time constant $T_{\phi}$ and, when not specified, a monolithic (SAND) solution strategy is adopted.
 %, while a normalized simulation time as normalized between $[0,1]$ where 0 corresponds to the initial time and 1 to the end simulation time $t_{end}$ at convergence.

\begin{table}[tbh]
\centering
\begin{tabular}{|c|c|c|c|c|c|c|c|c|c|c|}
\hline
Parameter & $a$ & $b$ & $c$ & $\Delta$ & $g$ & $E_A$ & $\nu_A$ & $\delta$ & $p$ \\ \hline
Unit & m & m & m & m & MPa & GPa & $-$ & $-$ & $-$  \\ \hline
Value & 2 & 1 & 1 & 0.2 & 500 & 10 & 0.25 & $10^{-3}$ & 10  \\ \hline
\end{tabular} 
\caption{Reference values of geometrical features, applied load, and material parameters of the 2D cantilever (Section \ref{sec:res2D}).
}
\label{tab:par1}
\end{table}

\begin{table}[tbh]
\centering
\begin{tabular}{|c|c|c|c|c|c|c|c|c|c|}
\hline
Parameter &  $\kappa_{\phi}$ & $\gamma_{\phi}$ & $\bar{v}$ & $\kappa_v$ & $\kappa_b$ &  $c_{AC}^{conv}$ \\ \hline
Unit   & MN/m & m  & $-$  & MPa & GPa & $\text{Pa s}$ \\ \hline
Value &  1 & 0.01 & 0.4 & 100 & $10^3$ &  $10^{-6}$ \\ \hline
\end{tabular} 
\caption{Reference values for simulation parameters adopted for the 2D cantilever (Section \ref{sec:res2D}).}
\label{tab:par2}
\end{table}

The design domain is discretized with $(120 \times 60)$ square elements along the two coordinate axes. As a result of a preliminary convergence analysis, the employed mesh density allows to obtain mesh-independent results for all cases under investigation.

Table \ref{tab:par2} reports the values of simulation parameters which are employed, if not object of a parametric investigation. 
It is noteworthy that adopted values respect the general rules-of-thumb traced in Section \ref{sec:par_settings}. In fact,  the  mesh size  is  $h_e =  0.016$ m and $\bar{e}^{el}_{,\phi} = 80$ MPa (with $g=500$ MPa). The target final volume fraction with the volume minimization functional hence results $v_{sol}^{tar} \approx 0.4$ with $\kappa_v=100$ MPa, that is comparable with the volume constraint value $\bar{v}=0.4$ assigned for the volume constraint functional.

\subsubsection{Staggered (NAND) versus monolithic (SAND) solution strategy: convergence behaviour}

The first set of analyses compares the convergence behaviour of a staggered (NAND) versus a monolithic (SAND) solution strategy. This comparison is made only with reference to the classical volume constraint approach, i.e., functional $\cL_{AC}^{vc}$.  Very  similar results are  however to be found  for the novel volume minimization approach, i.e., functional $\cL_{AC}^{vm}$,  as well. 

\begin{figure}[h!]
\centering
\begin{subfigure}[b]{0.99\textwidth}
\includegraphics[width=0.9\textwidth]{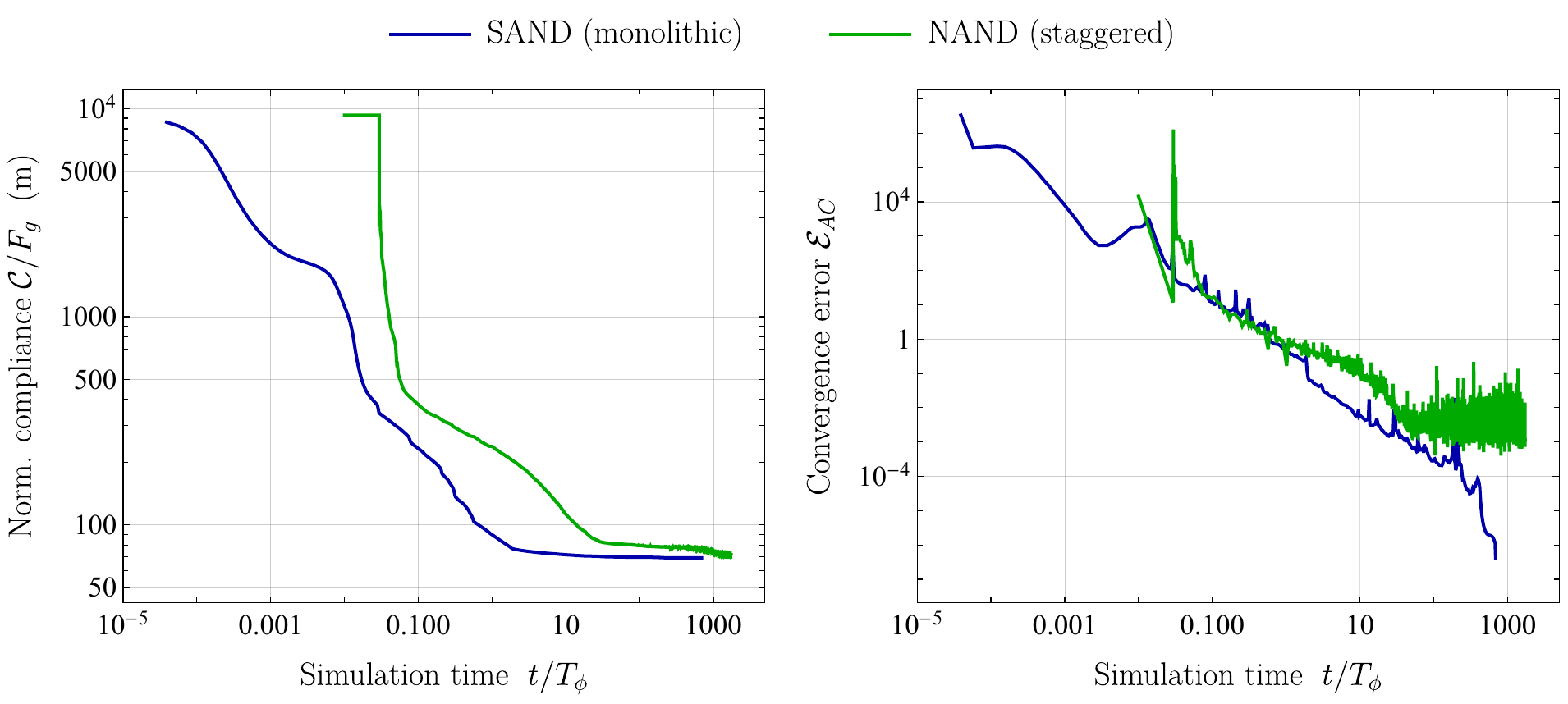}
\caption{Compliance $\mathcal{C}(\phi,{\bf u})$ (left) and Allen-Cahn error $\cE_{AC}$ (right).} \label{fig:conv_stagg_all_a}
 \end{subfigure}\vspace{2em}
 \begin{subfigure}[b]{0.99\textwidth}
\includegraphics[width=0.9\textwidth]{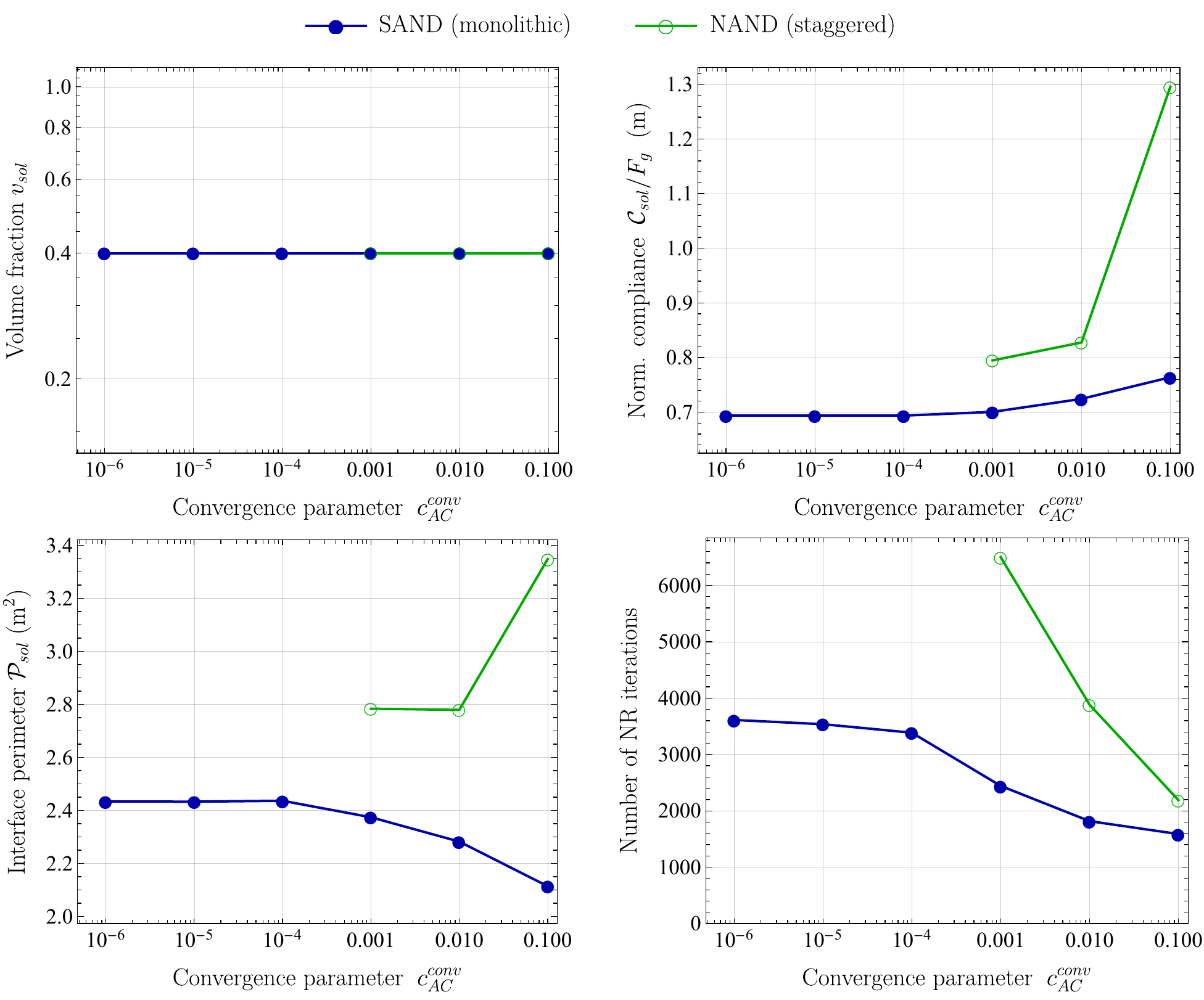}
\caption{Volume fraction $v_{sol}$ (top left), normalized  compliance $\cC_{sol}/F_g$  (top right), interface perimeter $\mathcal{P}_{sol}$  (bottom left), and total number of     
    Newton-Raphson iterations  (bottom right).} \label{fig:conv_stagg_all_b}
 \end{subfigure}
    \caption{
   Convergence behaviour of a staggered (SAND) and a monolithic solution (NAND) strategy: a) evolution of the computed solution versus the simulation time $t/T_{\phi}$;
    b) effect of the convergence parameter $c_{AC}^{conv}$ on the final computed solution. Results are shown only with reference to the classical volume constraint approach (i.e., functional $\cL_{AC}^{vc}$).
    }
\label{fig:conv_stagg_all}
\end{figure}

The evolution of the compliance $\cC$ and of  the Allen-Cahn error $\cE_{AC}$ versus the simulation time is shown in Fig. \ref{fig:conv_stagg_all_a}. It can be observed that the error of the staggered approach decreases at the beginning of the simulation but starts to oscillate when reaching $\cE_{AC}\approx 10^{-3}$, while the error of the monolithic strategy is substantially monotonic, down to  the  imposed value of $10^{-6}$.

\begin{figure}[h!]
\centering
\begin{subfigure}[tb]{0.99\textwidth}
\begin{subfigure}[tb]{0.32\textwidth}
\caption*{$c_{AC}^{conv}=10^{-3}$}
\includegraphics[width=0.99\textwidth]{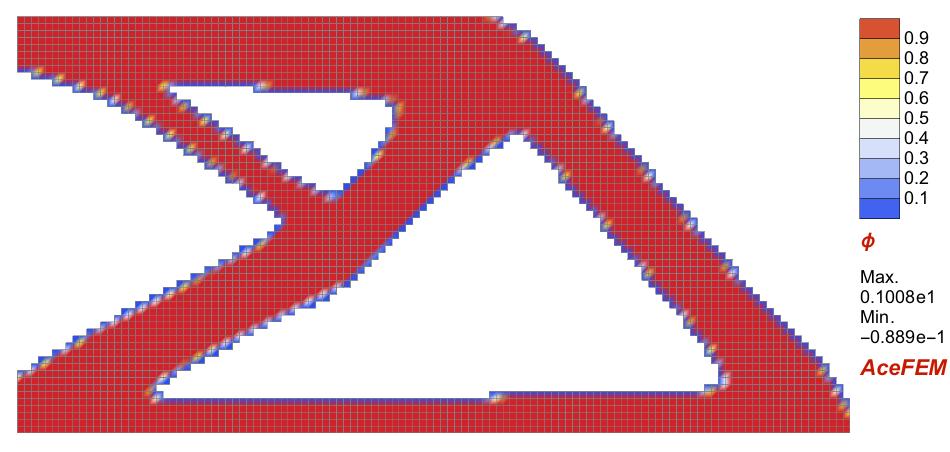}
 \end{subfigure}
\begin{subfigure}[tb]{0.32\textwidth}
\caption*{$c_{AC}^{conv}=10^{-2}$}
\includegraphics[width=0.99\textwidth]{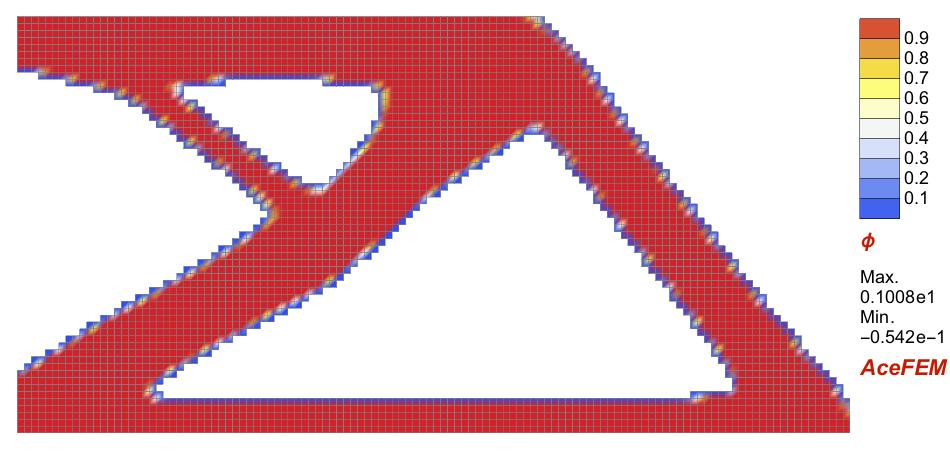}
 \end{subfigure}
 \begin{subfigure}[tb]{0.32\textwidth}
 \caption*{$c_{AC}^{conv}=10^{-1}$}
\includegraphics[width=0.99\textwidth]{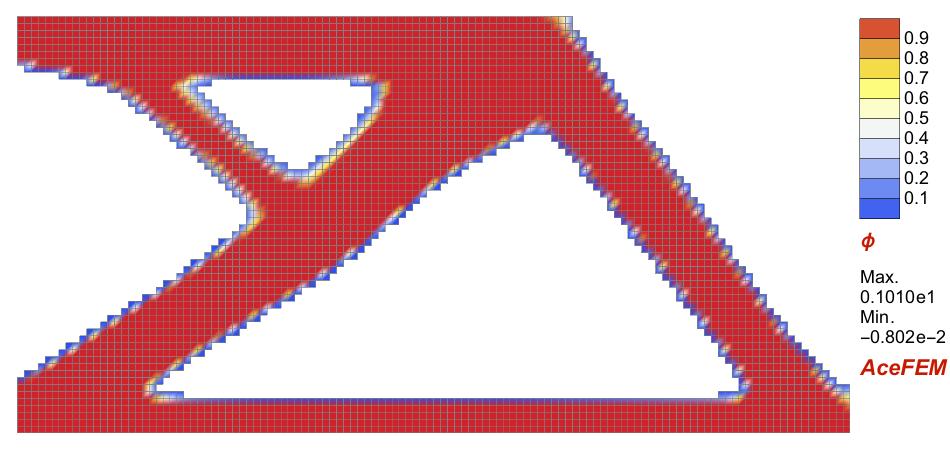}
 \end{subfigure}
 \caption{Monolithic approach (SAND)}
 \end{subfigure}

 \vspace{1em}
 \begin{subfigure}[tb]{0.99\textwidth}
 \begin{subfigure}[tb]{0.32\textwidth}
\includegraphics[width=0.99\textwidth]{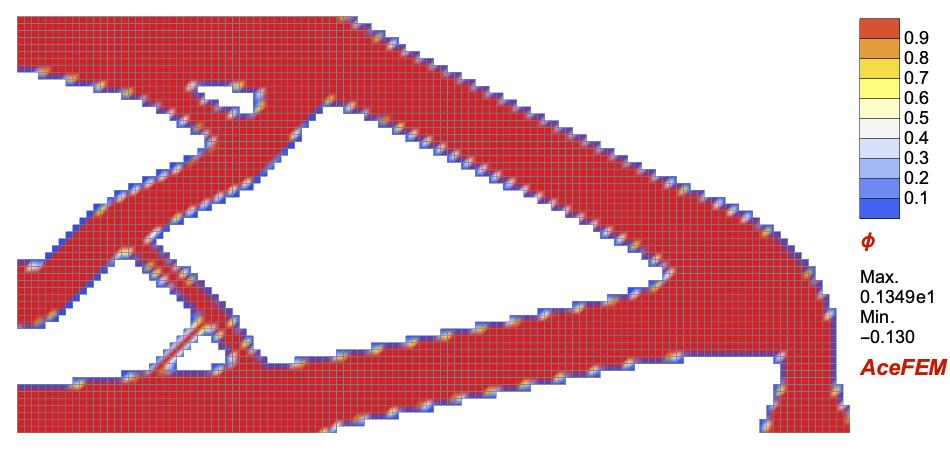}
 \end{subfigure}
\begin{subfigure}[tb]{0.32\textwidth}
\includegraphics[width=0.99\textwidth]{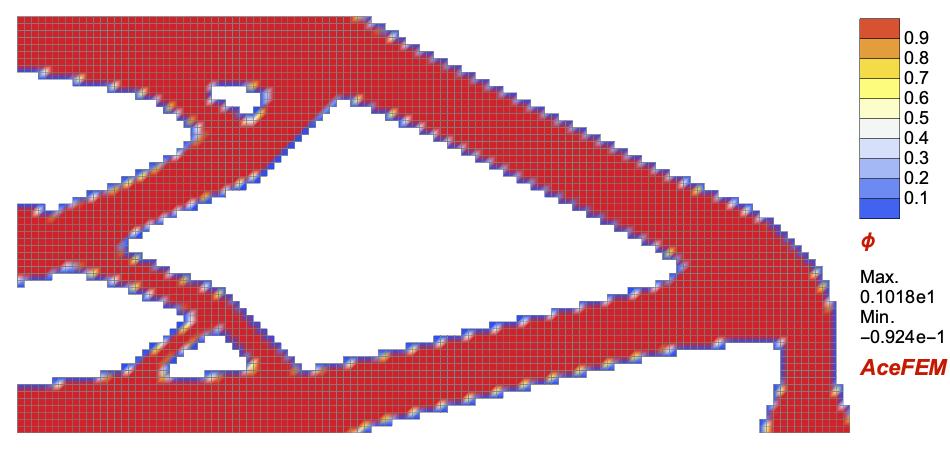}
 \end{subfigure}
 \begin{subfigure}[tb]{0.32\textwidth}
\includegraphics[width=0.99\textwidth]{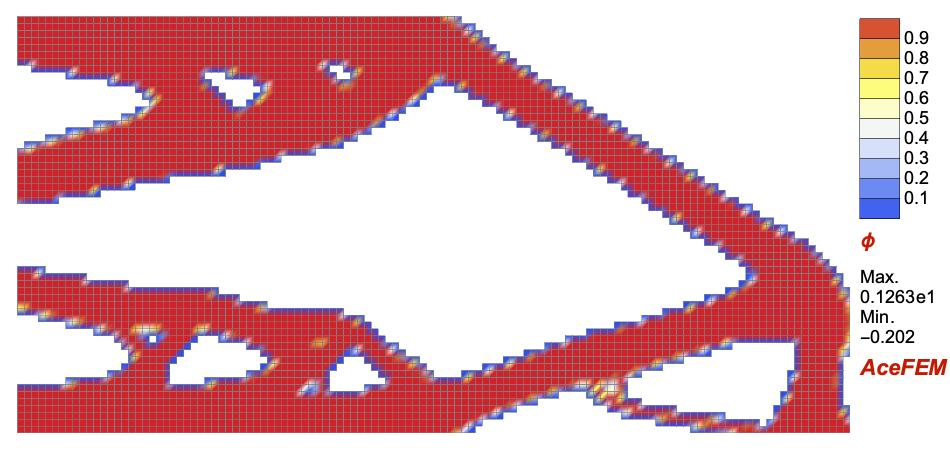}
 \end{subfigure}
 \caption{Staggered approach (NAND)}
 \end{subfigure}
    \caption{
  Material distribution phase field $\phi_{sol}$ obtained with different values of the convergence parameter $c_{AC}^{conv}$ for the monolithic (a) and the staggered (b) solution approaches. Results are shown only with reference to the classical volume constraint approach (i.e., functional $\cL_{AC}^{vc}$).
    }
\label{fig:conv_stagg_maps}
\end{figure}

The obtained final solution is then investigated for different values of the convergence parameter $c_{AC}^{conv}$,
value with respect to which convergence of the Allen-Cahn procedure is checked. 
Figure \ref{fig:conv_stagg_all_b} reports the volume fraction $v_{sol}$ and the  structural displacement $U_{sol}$ at the converged final solution, together with the final interface perimeter $\mathcal{P}_{sol}$. 
As previously noted, the staggered approach does not return a converged solution when enforcing $c_{AC}^{conv}\leq 10^{-3}$. The obtained solutions obviously coincide in terms of final volume fraction (since enforced by means of a global constraint), but not in terms of structural displacement and interface perimeter. A converged final  value  is not obtained with the staggered approach, while a converged solution is obtained with the monolithic one for $c_{AC}^{conv}\leq10^{-3}$, although $c_{AC}^{conv}=10^{-2}$ would be practically sufficient. This outcome is confirmed by the maps on the distribution of the phase field variable at different values of $c_{AC}^{conv}$, reported in Fig. \ref{fig:conv_stagg_maps}. 

Finally, Figure \ref{fig:conv_stagg_all_b} shows also the comparison on the total number of Newton-Raphson iterations, confirming the significant out-performance of the monolithic solution strategy (SAND) versus the staggered one (NAND). Therefore, all following analyses will adopt a monolithic (SAND) approach.

\subsubsection{Volume constraint versus volume minimization: convergence behaviour} \label{sec:convergence_res_form}

This section presents a comparison between the convergence behaviour of the volume constraint versus the volume minimization formulation. For both formulations, two different initial conditions are investigated, that is $\phi_0 = 1$ and $ \phi_0 = \bar{v}$ for $\mathcal{L}_{AC}^{vc}$, and $\phi_0 = 1$ and $ \phi_0 = 0.5$ for $\mathcal{L}_{AC}^{vm}$. 

\begin{figure}[h!]
\centering
\begin{subfigure}[b]{0.99\textwidth}
\includegraphics[width=0.9\textwidth]{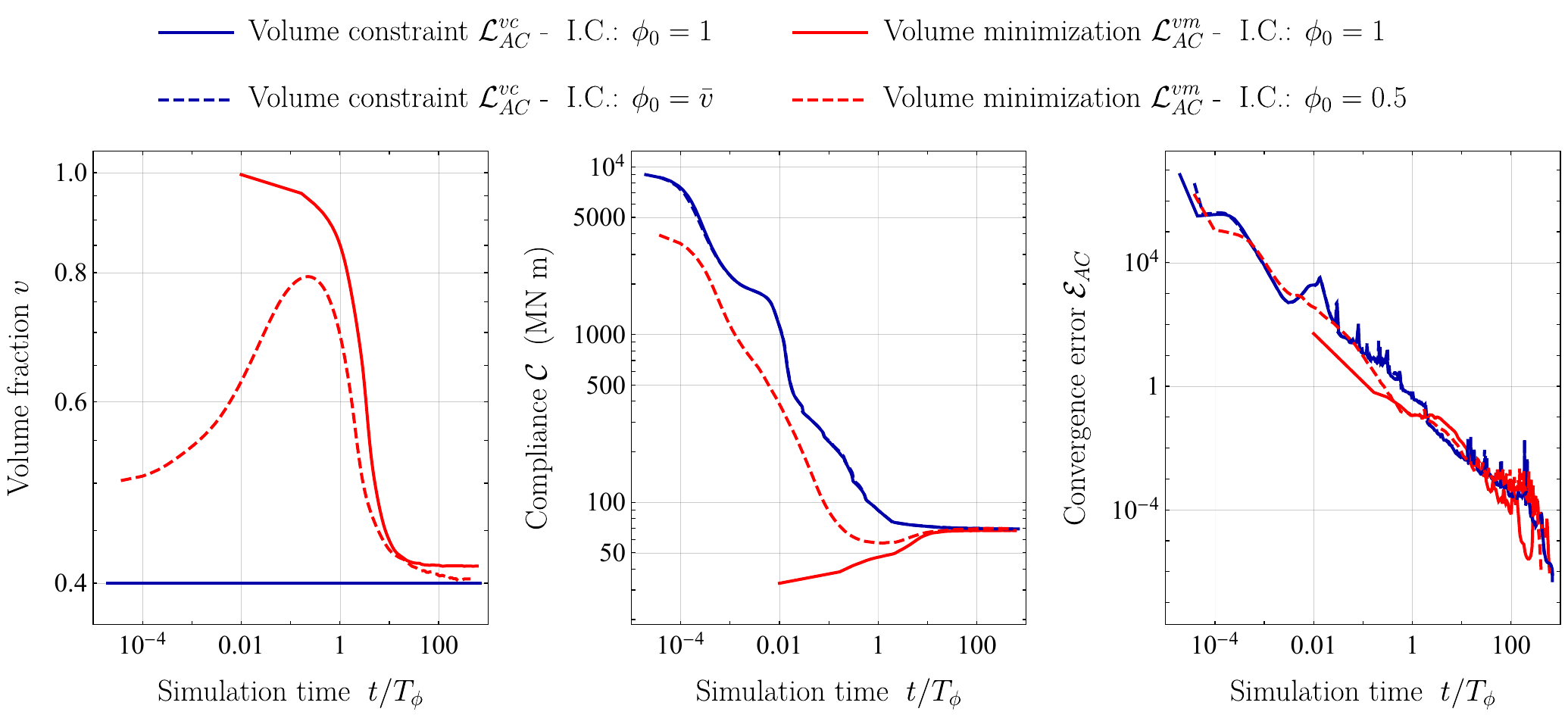}
\caption{Volume fraction $v(\phi)$ (left), compliance $\mathcal{C}(\phi,{\bf u})$ (center) and Allen-Cahn error $\cE_{AC}$ (right).} \label{fig:conv_all_a}
 \end{subfigure}\vspace{2em}
 \begin{subfigure}[b]{0.99\textwidth}
\includegraphics[width=0.9\textwidth]{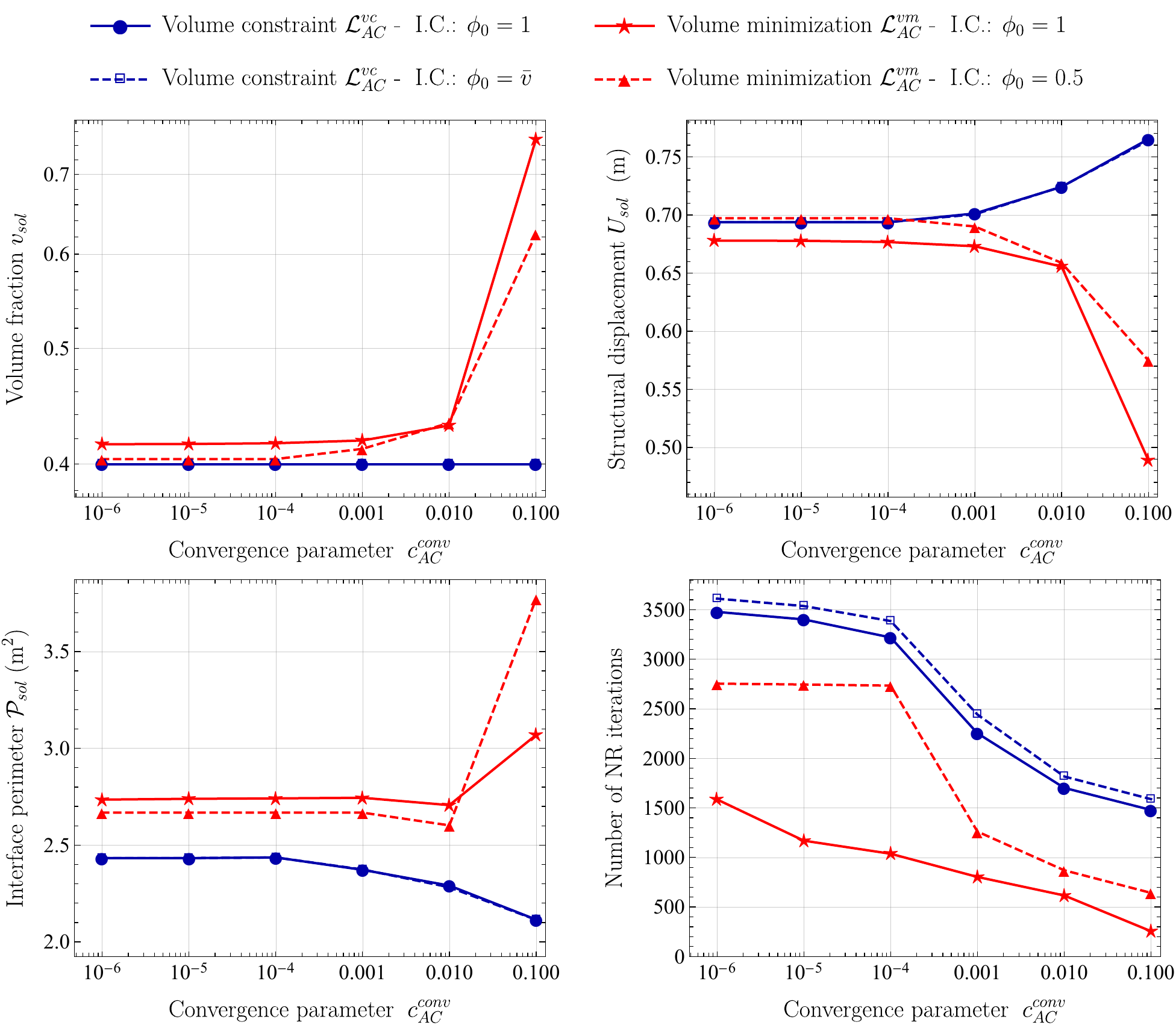}
\caption{Volume fraction $v_{sol}$ (top left), structural displacement $U_{sol}$  (top right), interface perimeter $\mathcal{P}_{sol}$  (bottom left), and total number of     
    Newton-Raphson iterations  (bottom right).} \label{fig:conv_all_b}
 \end{subfigure}
    \caption{
   Convergence behaviour of the two formulations with different initial conditions: a) 
   evolution of the computed solution versus the simulation time $t/T_{\phi}$;
    b) effect of the convergence parameter $c_{AC}^{conv}$ on the final computed solution.
    }
\label{fig:conv_all}
\end{figure}

Figure \ref{fig:conv_all_a} shows the evolution of the volume fraction $v$ and of the compliance   $\mathcal C$   along the simulation time. Contrarily to $\mathcal{L}_{AC}^{vc}$, the volume fraction $v$ with $\mathcal{L}_{AC}^{vm}$ evolves during the simulation, starting from the assigned initial condition. Structure compliance follows coherently. It is noteworthy that, for the chosen set of parameters, the final obtained solution between the two formulations is practically identical in terms of volume fraction and compliance, making the comparison consistent.

Figure \ref{fig:conv_all_a} reports also the Allen-Cahn convergence error $\mathcal{E}_{AC}$ along the simulation time. In all cases, results show an acceptable convergence behaviour, with a similar decrease rate of $\mathcal{E}_{AC}$  between the two formulations. However, the evolution of $\mathcal{E}_{AC}$ for the volume constraint functional $\mathcal{L}_{AC}^{vc}$ is significantly more oscillatory than the one obtained for the volume minimization one $\mathcal{L}_{AC}^{vm}$. 

\begin{figure}[h!]
\centering
\begin{subfigure}[tb]{0.99\textwidth}
\begin{subfigure}[tb]{0.32\textwidth}
\caption*{$c_{AC}^{conv}=10^{-6}$}
\includegraphics[width=0.99\textwidth]{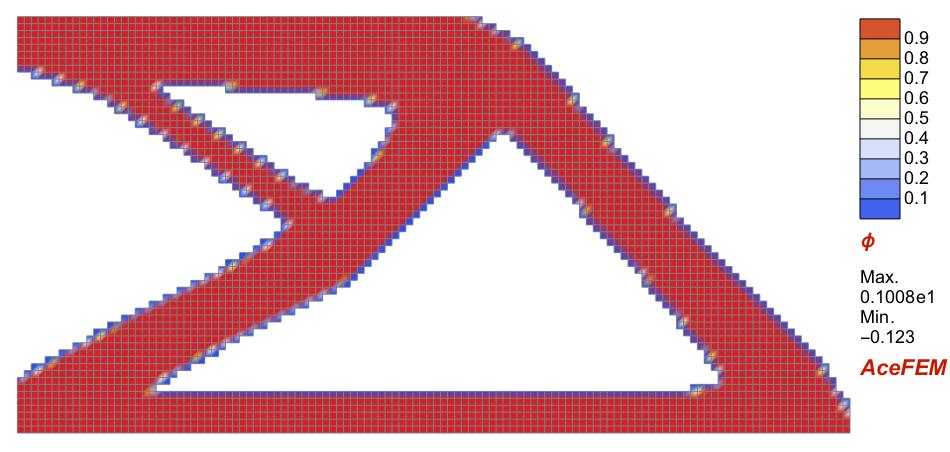}
 \end{subfigure}
\begin{subfigure}[tb]{0.32\textwidth}
\caption*{$c_{AC}^{conv}=10^{-4}$}
\includegraphics[width=0.99\textwidth]{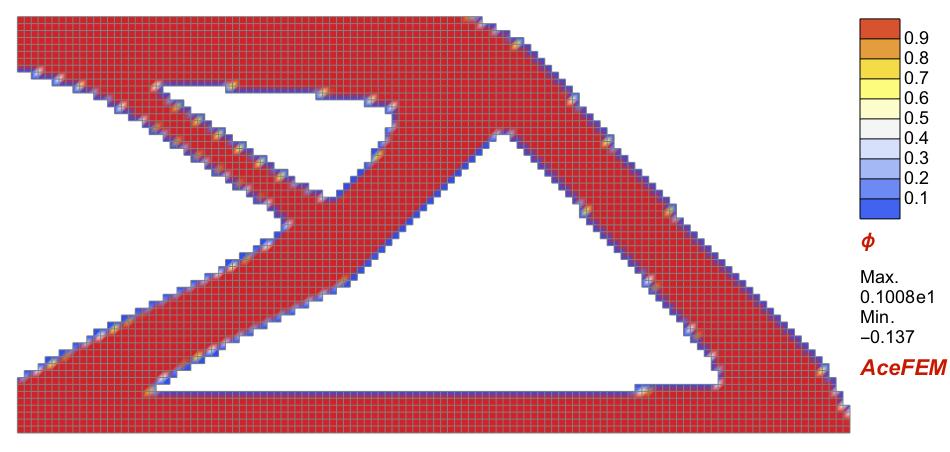}
 \end{subfigure}
 \begin{subfigure}[tb]{0.32\textwidth}
 \caption*{$c_{AC}^{conv}=10^{-2}$}
\includegraphics[width=0.99\textwidth]{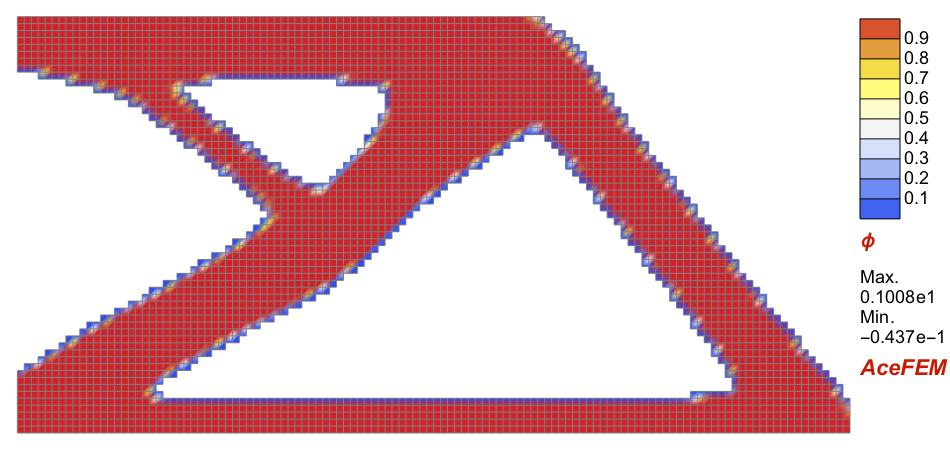}
 \end{subfigure}
 \caption{Volume constraint $\mathcal{L}_{AC}^{vc}$ $-$ $\text{I.C.}: \phi_0 = 1$ and $ \phi_0 = \bar{v}$}
 \end{subfigure}

 \vspace{1em}
 \begin{subfigure}[tb]{0.99\textwidth}
 \begin{subfigure}[tb]{0.32\textwidth}
\includegraphics[width=0.99\textwidth]{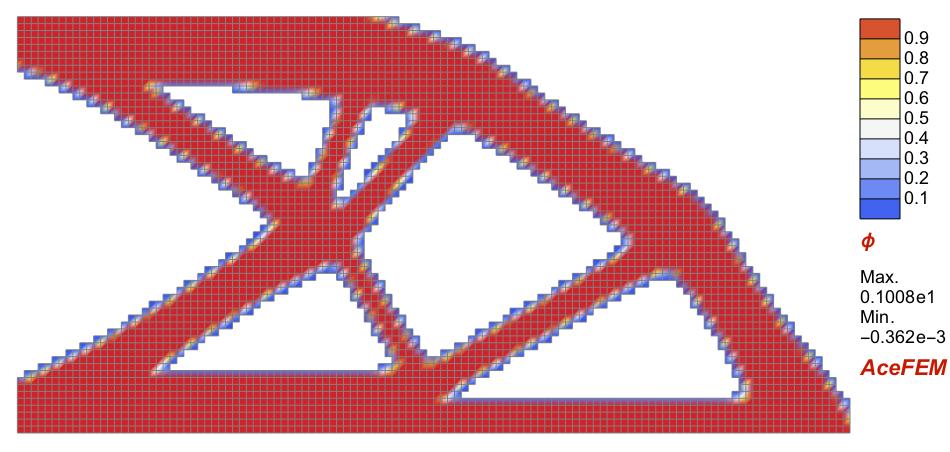}
 \end{subfigure}
\begin{subfigure}[tb]{0.32\textwidth}
\includegraphics[width=0.99\textwidth]{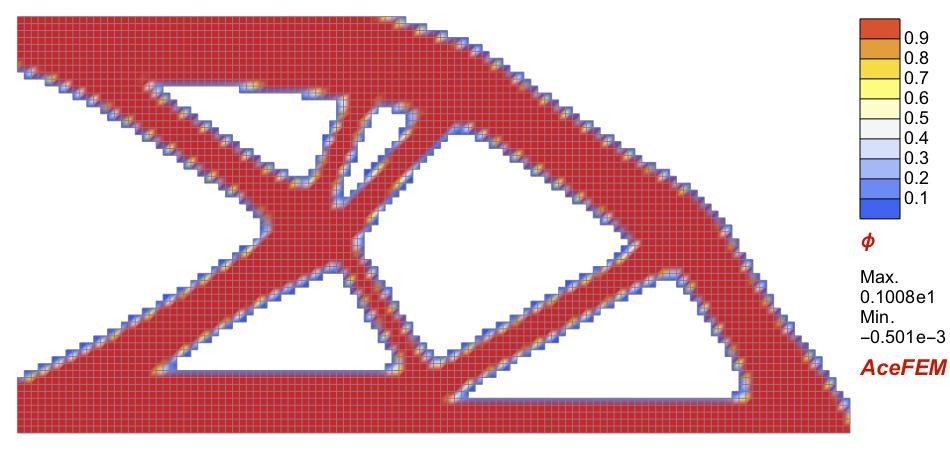}
 \end{subfigure}
 \begin{subfigure}[tb]{0.32\textwidth}
\includegraphics[width=0.99\textwidth]{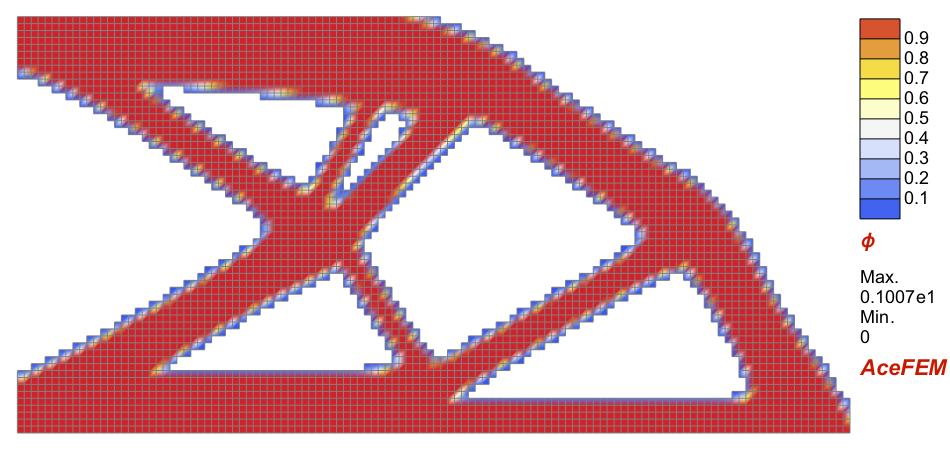}
 \end{subfigure}
 \caption{Volume minimization $\mathcal{L}_{AC}^{vm}$ $-$ $\text{I.C.}: \phi_0 = 1$}
 \end{subfigure}
 
  \vspace{1em}
 \begin{subfigure}[tb]{0.99\textwidth}
 \begin{subfigure}[tb]{0.32\textwidth}
\includegraphics[width=0.99\textwidth]{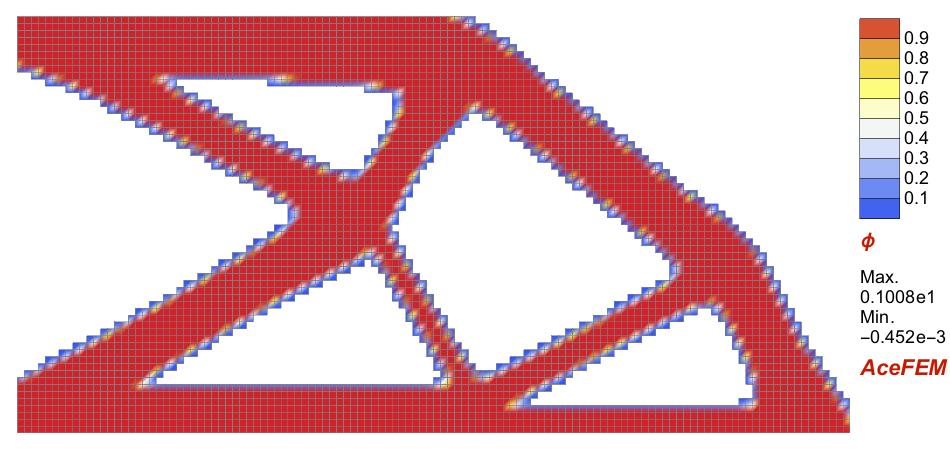}
 \end{subfigure}
\begin{subfigure}[tb]{0.32\textwidth}
\includegraphics[width=0.99\textwidth]{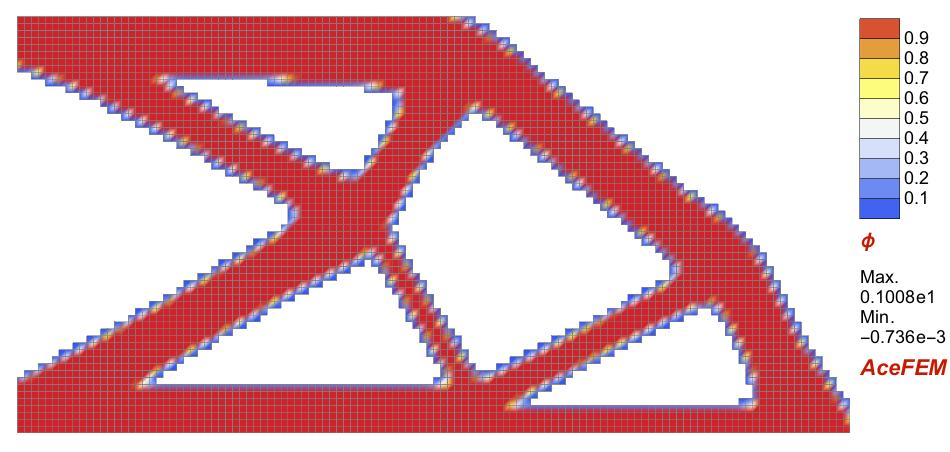}
 \end{subfigure}
 \begin{subfigure}[tb]{0.32\textwidth}
\includegraphics[width=0.99\textwidth]{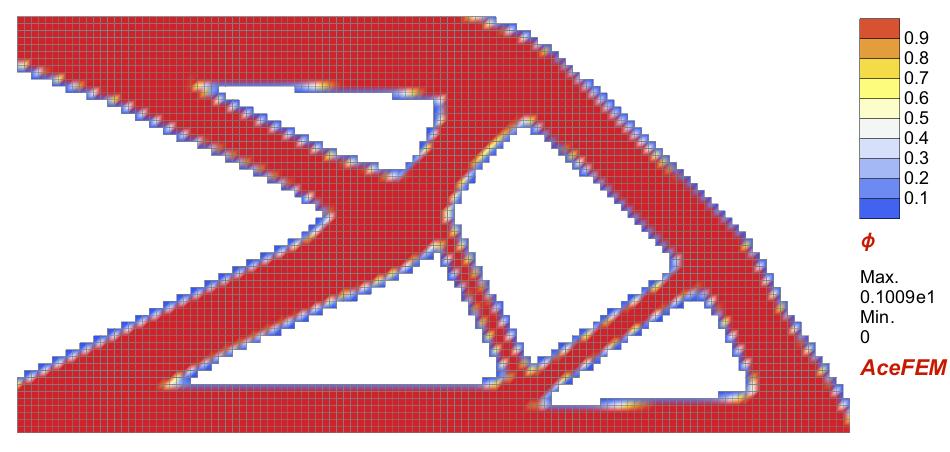}
 \end{subfigure}
 \caption{Volume minimization $\mathcal{L}_{AC}^{vm}$ $-$ $\text{I.C.}: \phi_0 = 0.5$}
 \end{subfigure}
    \caption{
  Material distribution phase field $\phi_{sol}$ obtained with different values of the convergence parameter $c_{AC}^{conv}$ and different initial conditions (I.C.) for: a) the volume constraint formulation (coinciding with $\phi_0 = \bar{v}$ and $\phi_0 = 1$); b) the volume minimization formulation with $\phi_0 = 1$; c) the volume minimization formulation with $\phi_0 = 0.5$. 
    }
\label{fig:conv_maps}
\end{figure}

The final  solution is then investigated for different values of the convergence parameter $c_{AC}^{conv}$ (with respect to which convergence of the Allen-Cahn procedure is checked). The final values of $v_{sol}$ and $U_{sol}$ are reported in Fig. \ref{fig:conv_all_b}, together with $\mathcal{P}_{sol}$. Below a given threshold (ca. $10^{-2}-10^{-4}$), both formulations show a robust and converged solution. This outcome is confirmed by the maps on the distribution of the phase field variable at different values of $c_{AC}^{conv}$, reported in Fig. \ref{fig:conv_maps}.

Moreover, Fig. \ref{fig:conv_all_b} reports the number of Newton-Raphson iterations required for the solution, clearly proving the computational out-performance of the volume minimization approach $\mathcal{L}_{AC}^{vm}$ with respect to the volume constraint one $\mathcal{L}_{AC}^{vc}$. 

The effect of the initial condition for the volume constraint functional $\mathcal{L}_{AC}^{vc}$ is negligible both in terms of final solution and numerical performances. This is due to the fact that, at the very first step, the global constraint changes the value of $\phi$ from $\phi_0$ to $\bar{v}$. On the other hand, for the volume minimization functional $\mathcal{L}_{AC}^{vm}$, the choice on $\phi_0$ affects the solution. In fact, only minor differences are obtained in terms of final obtained solution (see also Fig. \ref{fig:conv_maps}), but the convergence behaviour is significantly different, $\phi_0=1$ requiring a significantly lower number of Newton-Raphson iterations than $\phi_0=0.5$. This is also observable in Fig. \ref{fig:conv_all_a} by noting that the adaptive solution scheme with $\phi_0=0.5$  significantly reduces the reference time-step increment $\delta t_0 =10^{-2}$, while $\phi_0=1$ not. Accordingly, in what follows, the initial condition $\phi_0=1$ will be adopted for the volume minimization functional $\mathcal{L}_{AC}^{vm}$, while $\phi_0=\bar{v}$ for $\mathcal{L}_{AC}^{vc}$.

\subsubsection{Volume constraint versus volume minimization: amount of distributed material}

This section analyzes the two formulations by varying the amount of distributed material. This is achieved by varying $\bar{v}$ for the volume constraint functional $\cL_{AC}^{vc}$, and $\kappa_v$ for the volume minimization functional
 $\cL_{AC}^{vm}$. For the latter, the value of the perimeter stiffness parameter $\kappa_{\phi}$ varies with $\kappa_v$ according to the relationship $\kappa_{\phi}=\gamma_{\phi} \kappa_v$, as described in Section \ref{sec:par_settings}.

\begin{figure}[h!]
\centering
\includegraphics[width=0.9\textwidth]{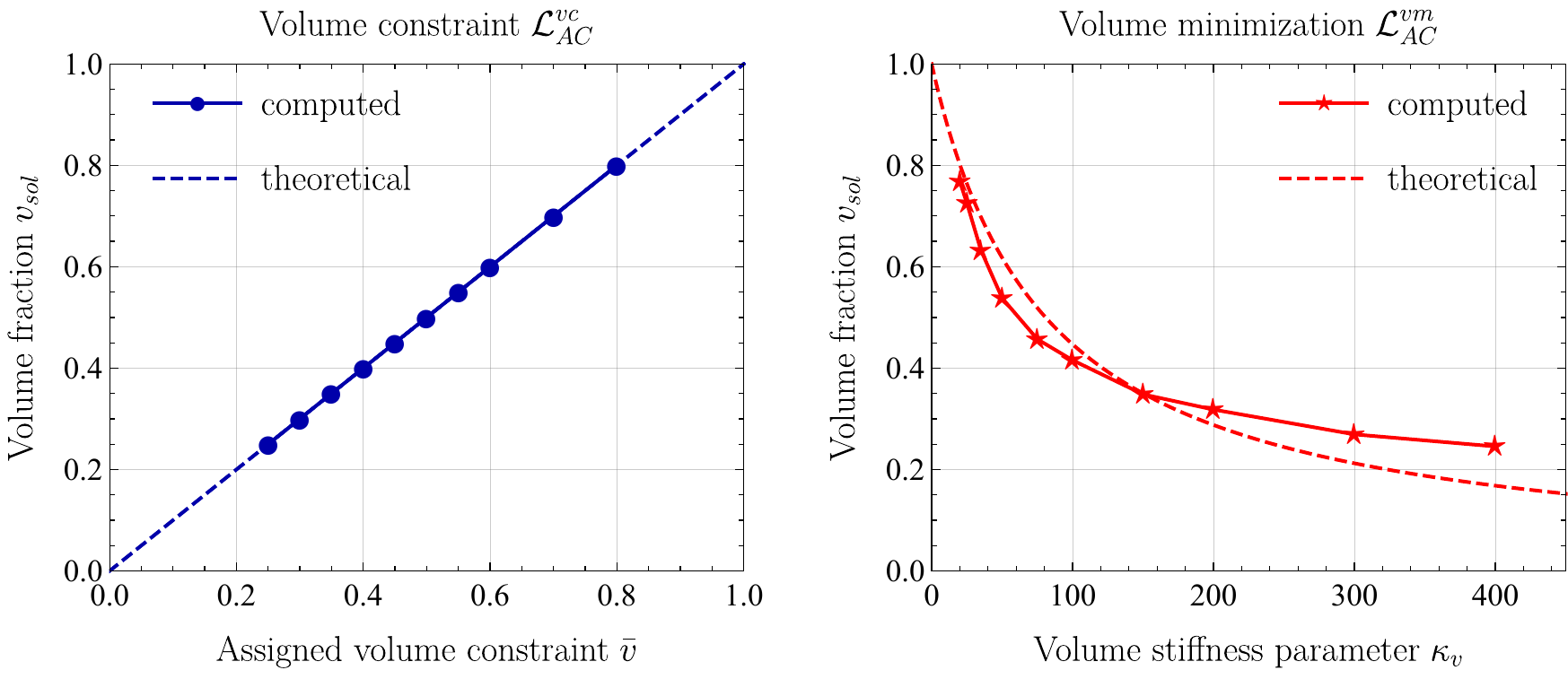}
\caption{Parametric study on the amount of distributed material with the volume constraint and the volume minimization functionals. Total volume fraction $v_{sol}$ obtained by varying the assigned volume constraint parameter $\bar{v}$ in $\mathcal{L}_{AC}^{vc}$ (left), and by varying the volume penalty parameter $\kappa_v$ (with $\kappa_{\phi}=\gamma_{\phi} \kappa_v$) in $\mathcal{L}_{AC}^{vm}$ (right). Computed values are compared with the estimated theoretical values $v_{sol}=\bar{v}$ for $\mathcal{L}_{AC}^{vc}$ and $v_{sol}=v_{sol}^{tar}(\kappa_v)$ in Eq. \eqref{eq:v_sol^est} for $\mathcal{L}_{AC}^{vm}$.}
\label{fig:mass_1}
\end{figure}

As shown in Fig. \ref{fig:mass_1}, for $\cL_{AC}^{vc}$, the assigned value $\bar{v}$ clearly corresponds to the obtained final volume fraction $v_{sol}$ for the entire range of the parametric analysis, verifying the correctness of the implementation of the global constraint. On the other hand, for $\cL_{AC}^{vc}$, the relationship between $\kappa_v$ and $v_{sol}$ is non-linear. Remarkably, the obtained trend is very-well captured by the theoretical target estimate $v_{sol}^{tar}(\kappa_v)$ provided in Eq. \eqref{eq:v_sol^est}. 

\begin{figure}[h!]
\centering
\includegraphics[width=0.9\textwidth]{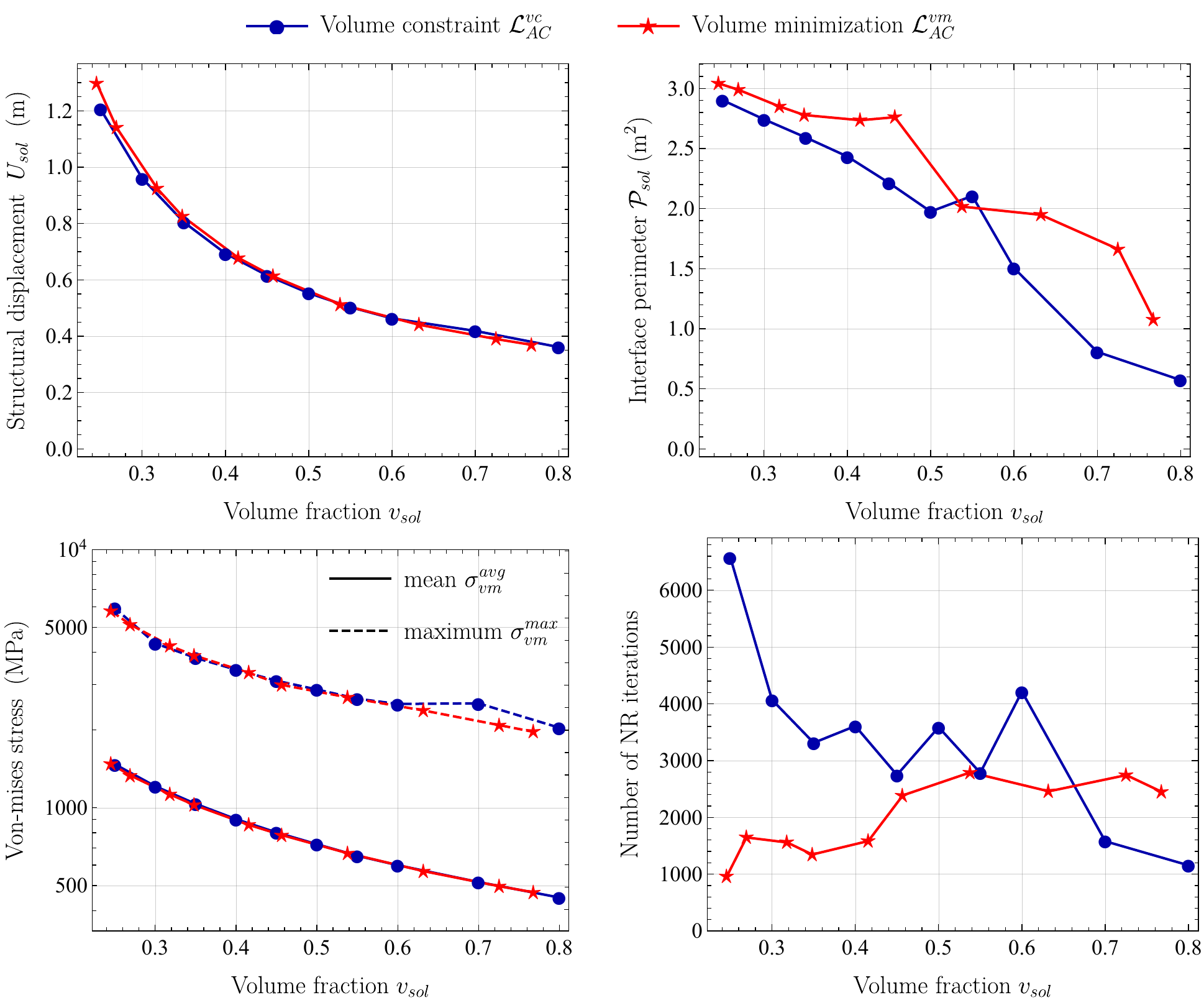}
\caption{Parametric study on the amount of distributed material with the volume constraint and the volume minimization functional. Structural displacement $U_{sol}$ (top left), interface perimeter $\mathcal{P}_{sol}$ (top right), maximum $\sigma_{vm}^{max}$  and average $\sigma_{vm}^{avg}$ Von-Mises stresses (bottom left), and total number of Newton-Raphson iterations (bottom right) versus the obtained total volume fraction $v_{sol}$.}
\label{fig:mass_2}
\end{figure}

\begin{figure}[h!]
\centering
\begin{subfigure}[tb]{0.99\textwidth}
\begin{subfigure}[tb]{0.32\textwidth}
\includegraphics[width=0.99\textwidth]{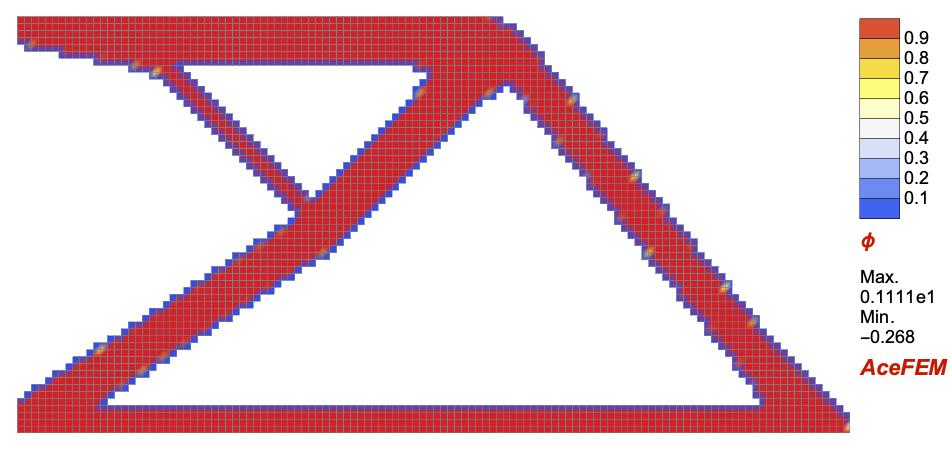}
\caption*{$\bar{v}=0.25$}
 \end{subfigure}
\begin{subfigure}[tb]{0.32\textwidth}
\includegraphics[width=0.99\textwidth]{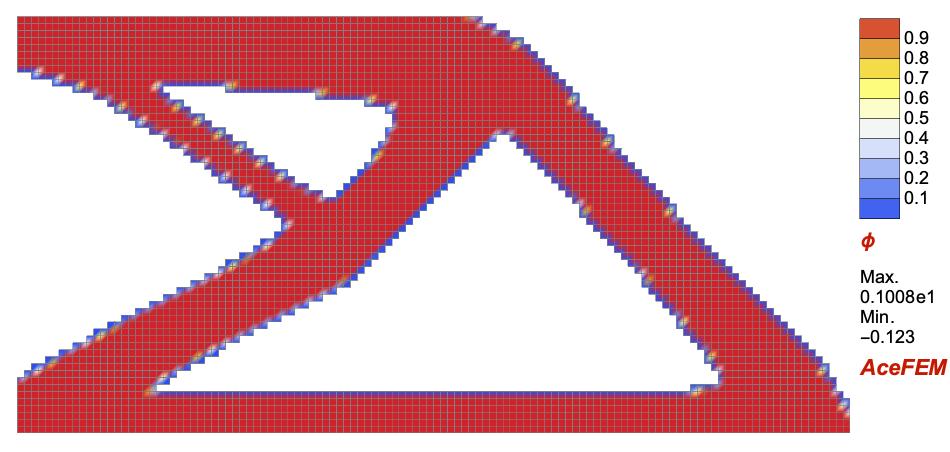}
\caption*{$\bar{v}=0.40$}
 \end{subfigure}
 \begin{subfigure}[tb]{0.32\textwidth}
\includegraphics[width=0.99\textwidth]{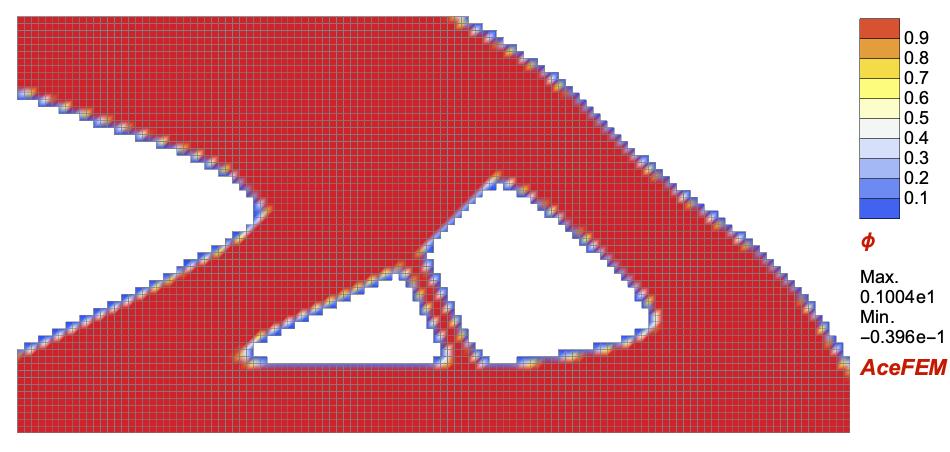}
\caption*{$\bar{v}=0.60$}
 \end{subfigure}
 \caption{Volume constraint $\mathcal{L}_{AC}^{vc}$}
 \end{subfigure}

 \vspace{1em}
 \begin{subfigure}[tb]{0.99\textwidth}
\begin{subfigure}[tb]{0.32\textwidth}
\includegraphics[width=0.99\textwidth]{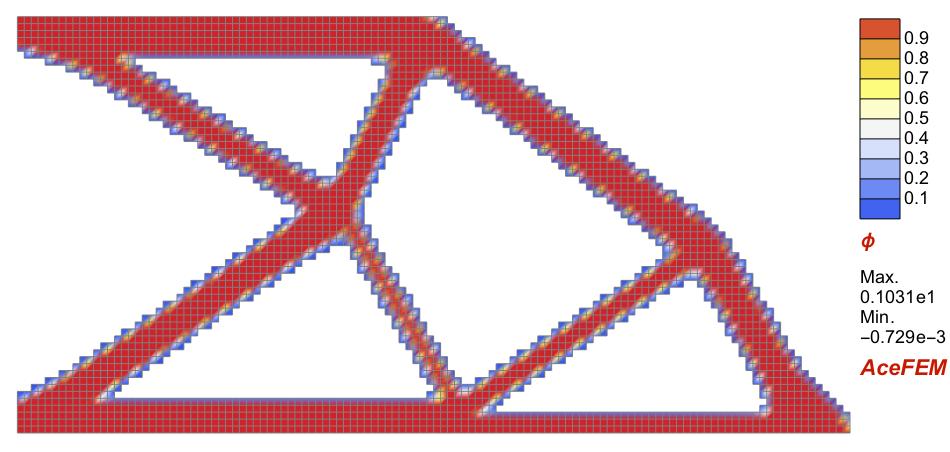}
\caption*{$\kappa_v = 300$ MPa, $v_{sol}\approx 0.25$}
 \end{subfigure}
\begin{subfigure}[tb]{0.32\textwidth}
\includegraphics[width=0.99\textwidth]{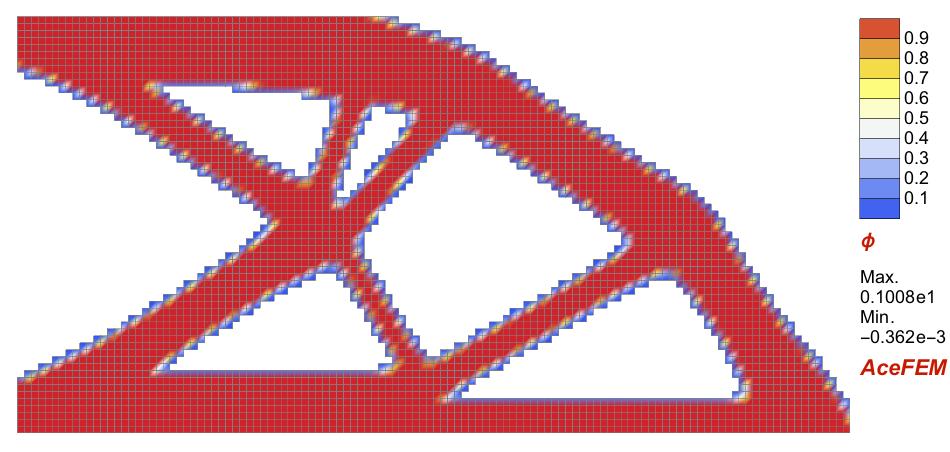}
\caption*{$\kappa_v = 100$ MPa, $v_{sol}\approx 0.4$}
 \end{subfigure}
 \begin{subfigure}[tb]{0.32\textwidth}
\includegraphics[width=0.99\textwidth]{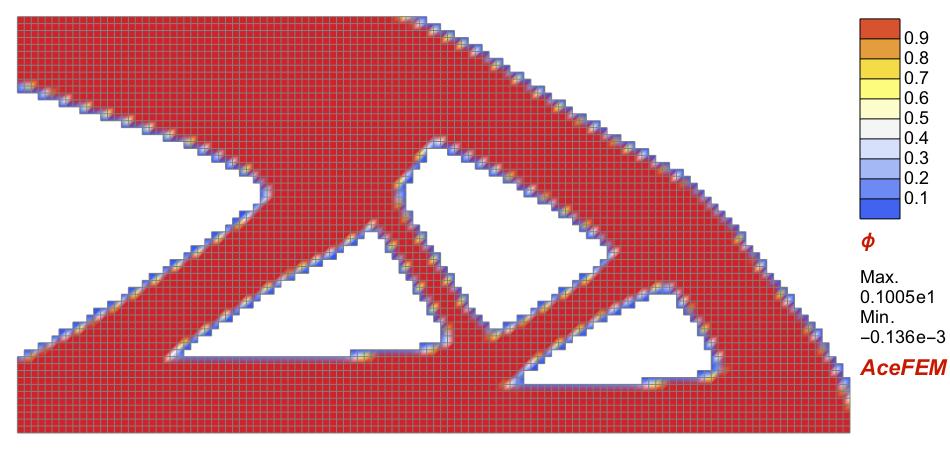}
\caption*{$\kappa_v = 50$ MPa, $v_{sol}\approx 0.6$}
 \end{subfigure}
 \caption{Volume minimization $\mathcal{L}_{AC}^{vm}$}
 \end{subfigure}
    \caption{
  Material distribution phase field $\phi_{sol}$ obtained with different values of obtained by varying the assigned volume constraint parameter $\bar{v}$ in $\mathcal{L}_{AC}^{vc}$ (a), and by varying the volume penalty parameter $\kappa_v$ (with $\kappa_{\phi}=\gamma_{\phi} \kappa_v$) in $\mathcal{L}_{AC}^{vm}$ (b).
    }
\label{fig:mass_maps}
\end{figure}

Figure \ref{fig:mass_2} shows that, for a given value of final volume fraction $v_{sol}$, the solutions obtained from the two formulations are practically identical in terms of structural displacement $U_{sol}$, as well as maximum $\sigma_{vm}^{max}$  and average $\sigma_{vm}^{avg}$ Von-Mises stresses. The effective distribution of the phase-field variable is also very similar, although differences are observable in Fig. \ref{fig:mass_maps} and again in Fig. \ref{fig:mass_2} as regards the interface perimeter $\mathcal{P}_{sol}$. 

Finally, the analysis of the performance of the two formulations (in terms of Newton-Raphson iterations, see Fig. \ref{fig:mass_2}) show that the volume minimization functional $\cL_{AC}^{vm}$ is significantly more efficient ($>50\%$ iteration saving) than $\cL_{AC}^{vc}$ for the more challenging cases cases where the final volume fraction is low, that is $v_{sol}<0.5$, cases which are also more interesting from the engineering viewpoint.

\subsubsection{Volume constraint versus volume minimization: effect of load variations}

The last comparison between the two formulations addresses the effect of the applied load magnitude $g$ on the final solution. In this case, the volume constraint $\bar{v}$ (for $\cL_{AC}^{vc}$) and the volume penalty $\kappa_v$ (for $\cL_{AC}^{vm}$) are held constant, as given in Table \ref{tab:par2}. 

\begin{figure}[h!]
\centering
\includegraphics[width=0.9\textwidth]{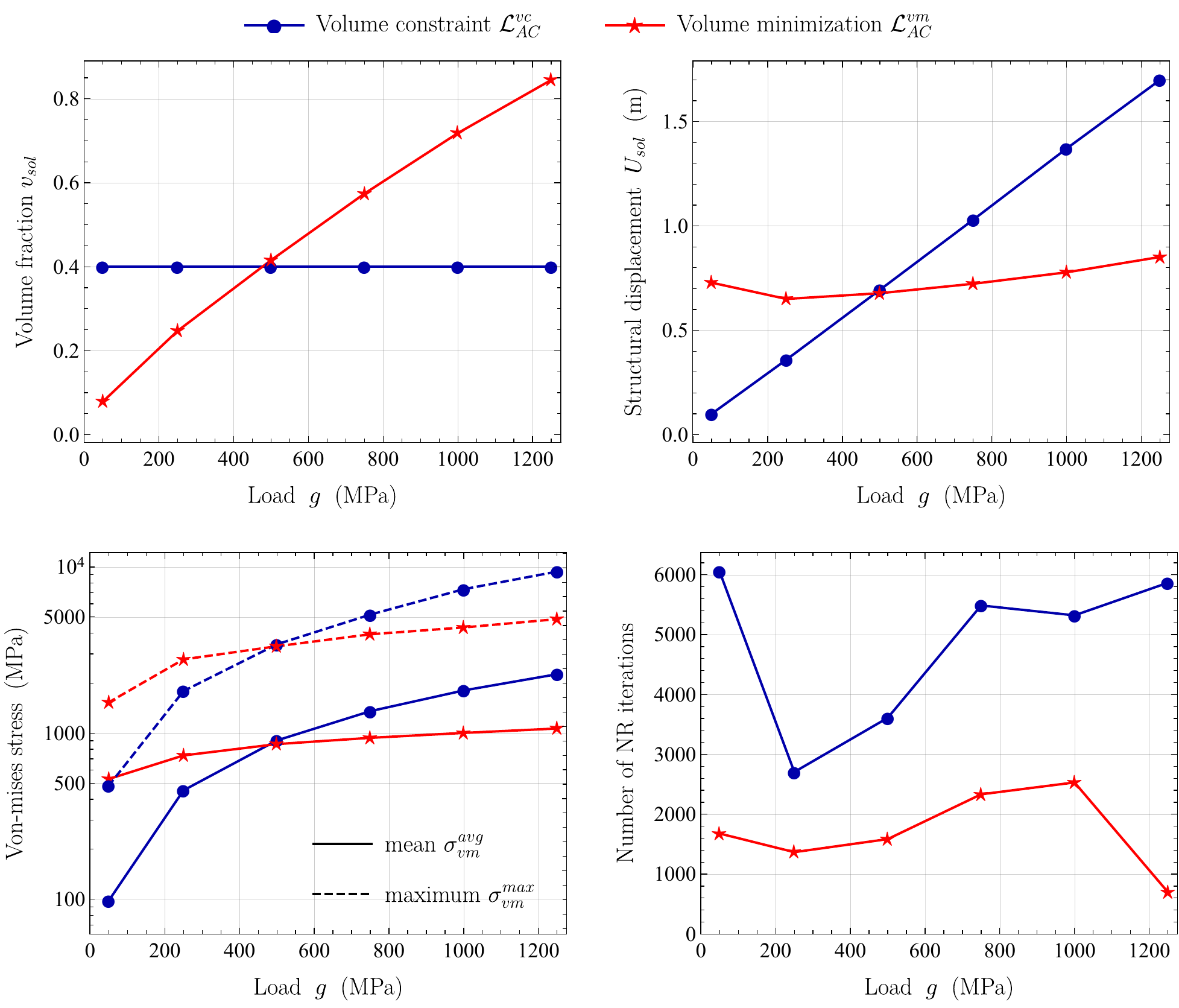}
\caption{Parametric study on the applied load with the volume constraint ($\bar{v}=0.4$) and the volume minimization  ($\kappa_v=100$ MPa) functionals. Total volume fraction $v_{sol}$ (top left), structural displacement $U_{sol}$ (top right), maximum $\sigma_{vm}^{max}$ and average $\sigma_{vm}^{avg}$ Von-Mises stresses (bottom left), and total number of Newton-Raphson iterations (bottom right) versus the applied load magnitude $g$.}
\label{fig:load}
\end{figure}

\begin{figure}[h!]
\centering
\begin{subfigure}[tb]{0.99\textwidth}
\begin{subfigure}[tb]{0.32\textwidth}
\caption*{$g=250$ MPa}
\includegraphics[width=0.99\textwidth]{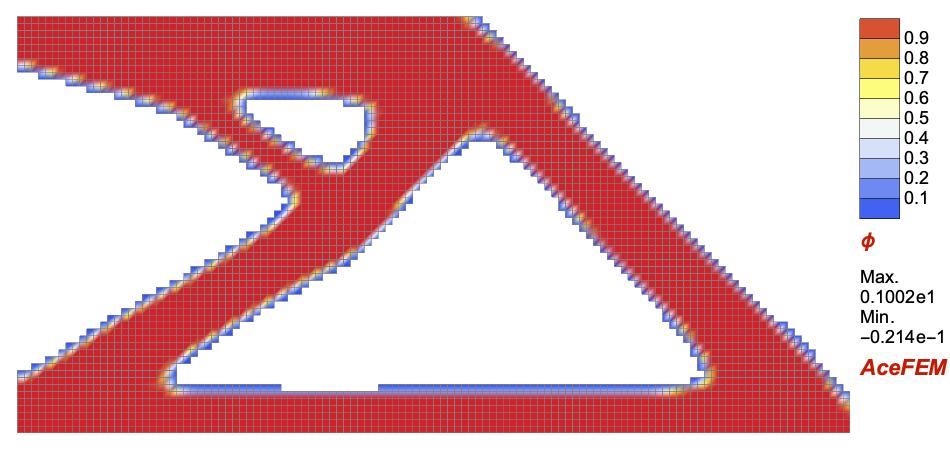}
 \end{subfigure}
\begin{subfigure}[tb]{0.32\textwidth}
\caption*{$g=500$ MPa}
\includegraphics[width=0.99\textwidth]{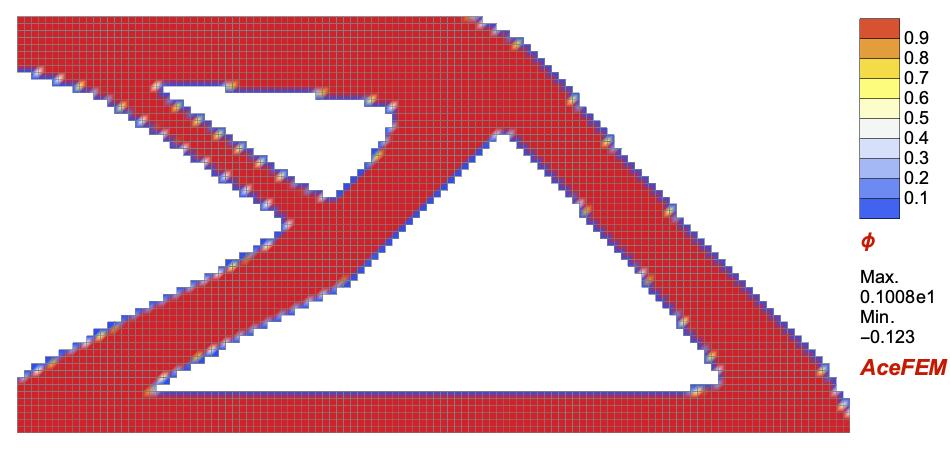}
 \end{subfigure}
 \begin{subfigure}[tb]{0.32\textwidth}
 \caption*{$g=750$ MPa}
\includegraphics[width=0.99\textwidth]{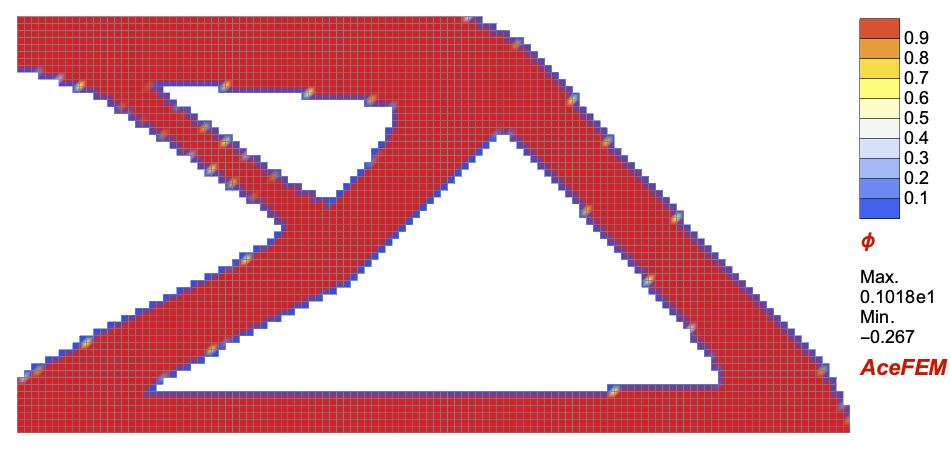}
 \end{subfigure}
 \caption{Volume constraint $\mathcal{L}_{AC}^{vc}$}
 \end{subfigure}

 \vspace{1em}
 \begin{subfigure}[tb]{0.99\textwidth}
\begin{subfigure}[tb]{0.32\textwidth}
\includegraphics[width=0.99\textwidth]{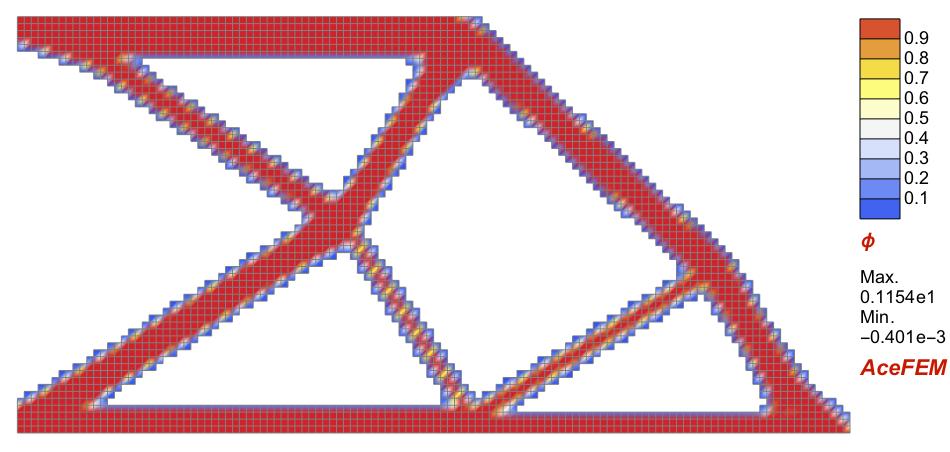}
 \end{subfigure}
\begin{subfigure}[tb]{0.32\textwidth}
\includegraphics[width=0.99\textwidth]{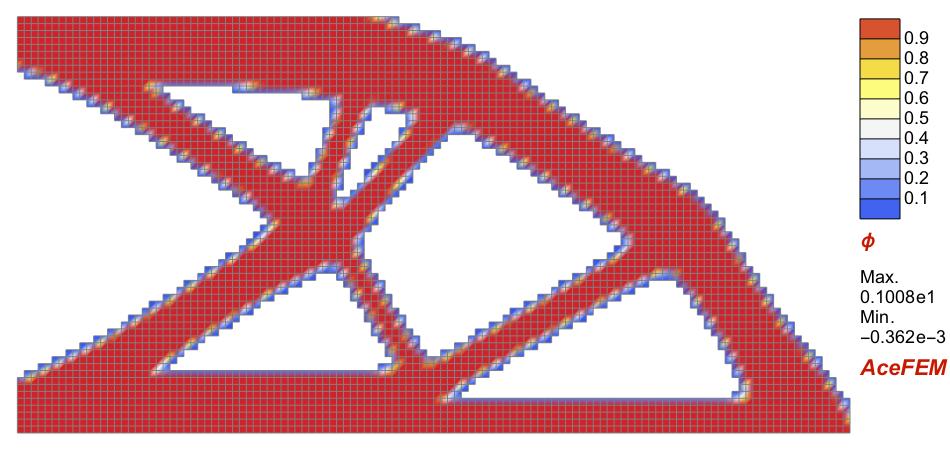}
 \end{subfigure}
 \begin{subfigure}[tb]{0.32\textwidth}
\includegraphics[width=0.99\textwidth]{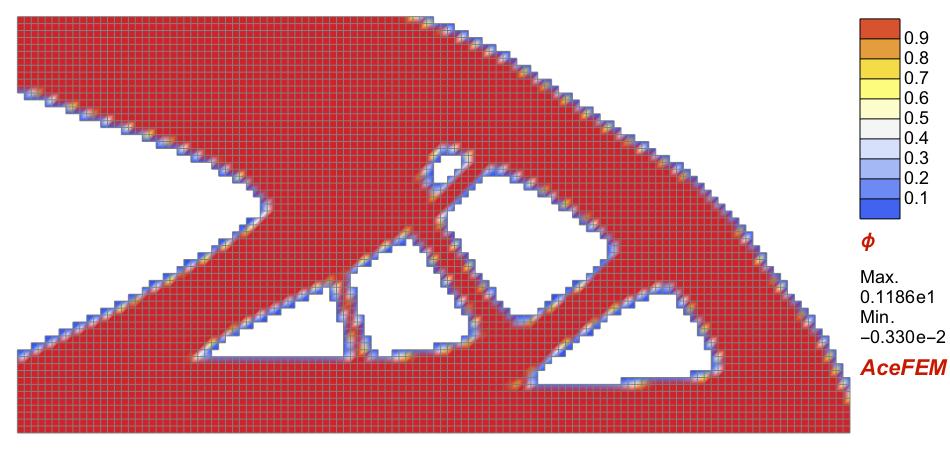}
 \end{subfigure}
 \caption{Volume minimization $\mathcal{L}_{AC}^{vm}$}
 \end{subfigure}
    \caption{
  Material distribution phase field $\phi_{sol}$ obtained with different values of the applied load magnitude $g$ for: a) the volume constraint formulation with $\bar{v}=0.4$; b) the volume minimization formulation with $\kappa_v=100$ MPa.
    }
\label{fig:load_maps}
\end{figure}

Figure \ref{fig:load} clearly show the difference on the rationale upon which the two functionals are built. In fact,  the  load magnitude $g$ highly affects the final obtained structural displacement $U_{sol}$ when employing a volume constrained principle (i.e., $\cL_{AC}^{vc}$) because the final volume fraction $v_{sol}$ is assigned \emph{a priori}. This is accompanied by high variations in stresses, while the final design is practically independent from the load (see Fig. \ref{fig:load_maps}). On the other hand,   the   load magnitude $g$ highly affects the final volume fraction $v_{sol}$ obtained when employing a volume minimization principle (i.e., $\cL_{AC}^{vm}$). At different load levels, final designs associated with similar values of structural displacement $U_{sol}$ and stresses  $\sigma_{vm}^{max}$ and $\sigma_{vm}^{avg}$ are obtained. As shown in Fig. \ref{fig:load_maps}, the obtained design with $\cL_{AC}^{vm}$ is now highly affected by the applied load value $g$.

\clearpage

\clearpage

\clearpage

\clearpage 

\clearpage

\subsection{3D simulations} \label{sec:res3D}

In this section, the novel volume minimization functional $\cL_{AC}^{vm}$ is tested on two 3D problems. In particular, the addressed examples are (see Fig. \ref{fig:test1_geometry}):

\begin{itemize}

\item A cantilever beam, representing the 3D counterpart of the case study presented in Section \ref{sec:res2D}. The domain is now discretized by means of  $80 \times 40 \times 40$  hexahedral elements, resulting in an average element size  $h_e = 0.025$  m. The geometrical and material features, as well as the applied load magnitude, correspond to the values listed in Table \ref{tab:par1}, while values of simulation parameters to the ones in Table \ref{tab:par2} apart from $\gamma_{\phi}=0.02$ m (due to a slight increase in the dimension of the elements in the mesh), and $c_{AC}^{conv}=10^{-2}$ (fixed on the basis of the convergence analysis in Section \ref{sec:convergence_res_form});

\item A bridge, as the one investigated in \cite{Zegard2016}. The design domain is represented by a rectangular cuboid  $\ell_x \times \ell_y \times \ell_z$. The domain is fixed on the bottom plane at $Z=0$ along strips of length $l_{bc}$ along $X$. Moreover, a passive solid slab on the top surface at $Z=\ell_z$ of height $h_s$ m is considered, on top of which a distributed load $g$ is applied in the $Z$-direction. Material properties correspond to the one of the previous example, being reported in Table~\ref{tab:par1}. Values of geometrical and simulation parameters, as well as load magnitude, are given in Table \ref{tab:par_bridge}. The domain is discretized by means of $150 \times 30 \times 30$ hexahedral elements, resulting in an average element size $h_e\approx 0.3$ m.  

\end{itemize}

\begin{table}[tbh]
\centering
\begin{tabular}{|c|c|c|c|c|c|c|c|c|c|c|c|c|c|c|c|c|c|c|c|}
\hline
Parameter & $\ell_x$ & $\ell_y$ & $\ell_z$ & $\ell_{bc}$ & $h_s$ &$g$ & $\kappa_{\phi}$ & $\gamma_{\phi}$ & $\kappa_v$  & $c_{AC}^{conv}$ \\ \hline
Unit 	& m & m & m & m & m & MPa & MN/m& m  &  MPa & $-$\\ \hline
Value & 44 & 8.8 & 8.8 & 1.76 & 0.2 & 150 &  200 & 0.2 & 1000 &  $10^{-2}$ \\ \hline
\end{tabular} 
\caption{Values of geometrical and simulation parameters employed in the 3D bridge case study. Only values different from the ones in Tables \ref{tab:par1} and \ref{tab:par2} are reported.}
\label{tab:par_bridge}
\end{table}

In particular, the 3D cantilever  has  been addressed to show that the proposed formulation can be straightforwardly implemented in a 3D computational framework, analyzing differences  obtained from a 2D to a 3D setting for the same case study. The bridge is instead defined such to have significant differences in the physical dimension of the problem, in order to prove the robustness of the parameter settings guidelines traced in Section \ref{sec:par_settings}, with particular reference to $\kappa_v$, $\gamma_{\phi}$, and $\kappa_{\phi}$.

\begin{figure}[tb]
\centering
\includegraphics[width=0.85\textwidth]{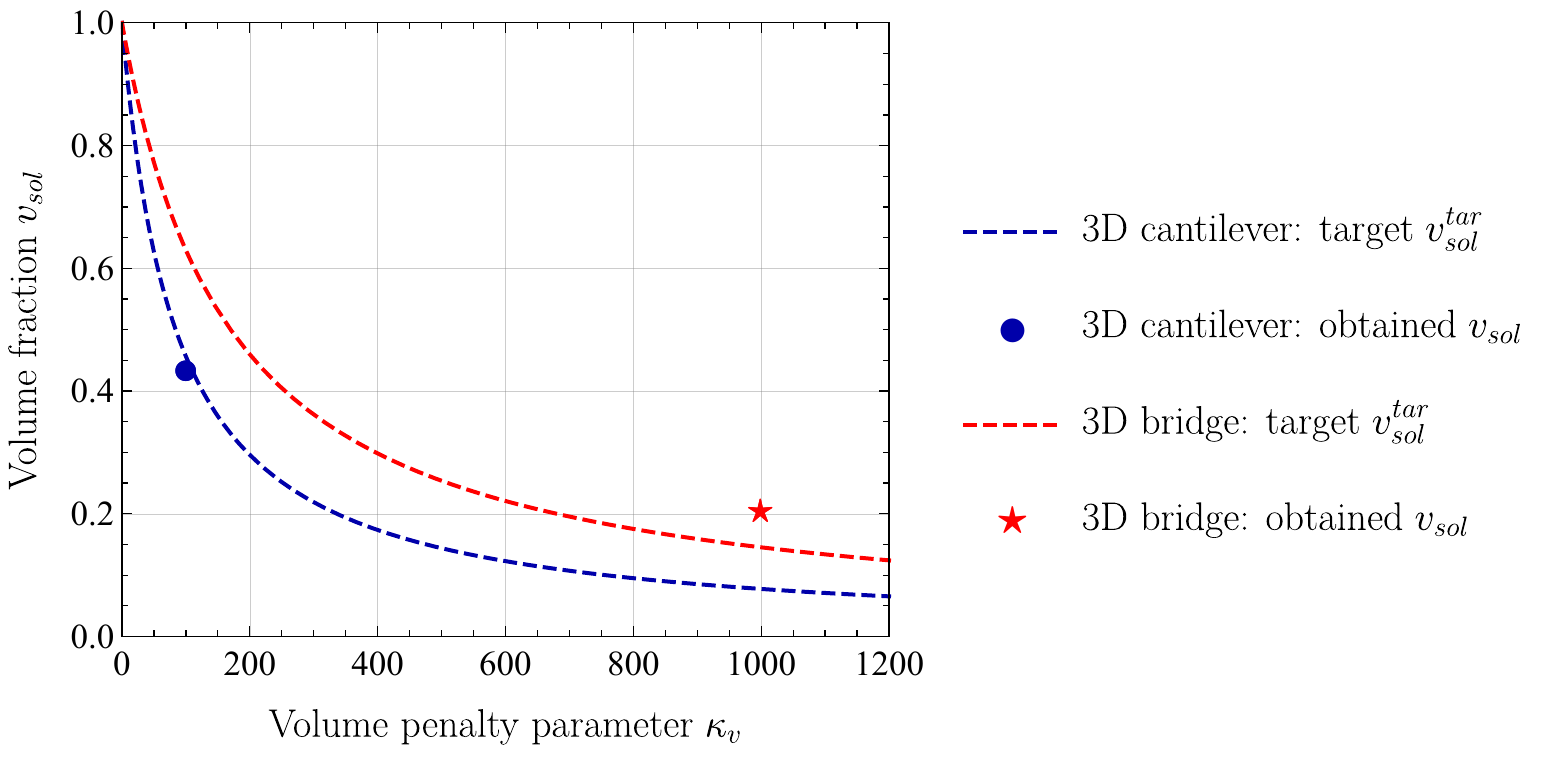}
        \caption{
        3D applications: estimated target volume fractions  $v_{sol}^{tar}$ versus obtained volume fraction $v_{sol}$. }
        \label{fig:mass3D}
\end{figure}

Since the volume minimization functional is adopted, the setting of the volume penalty parameter plays a fundamental role. For the cantilever,  we have  $\bar{e}^{el}_{\phi}=84$ MPa and hence $v_{sol}^{tar} \approx 0.4$ with $\kappa_v=100$ MPa; for the bridge,  we have   $\bar{e}^{el}_{\phi}=170$ MPa and hence $v_{sol}^{tar} \approx 0.2$ with $\kappa_v=1000$ MPa. Figure \ref{fig:mass3D} shows the comparison between the estimated and final obtained volume fractions, confirming the effectiveness of the proposed procedure for obtaining the target solution.

\begin{figure}[h!]
\centering
\includegraphics[width=0.99\textwidth]{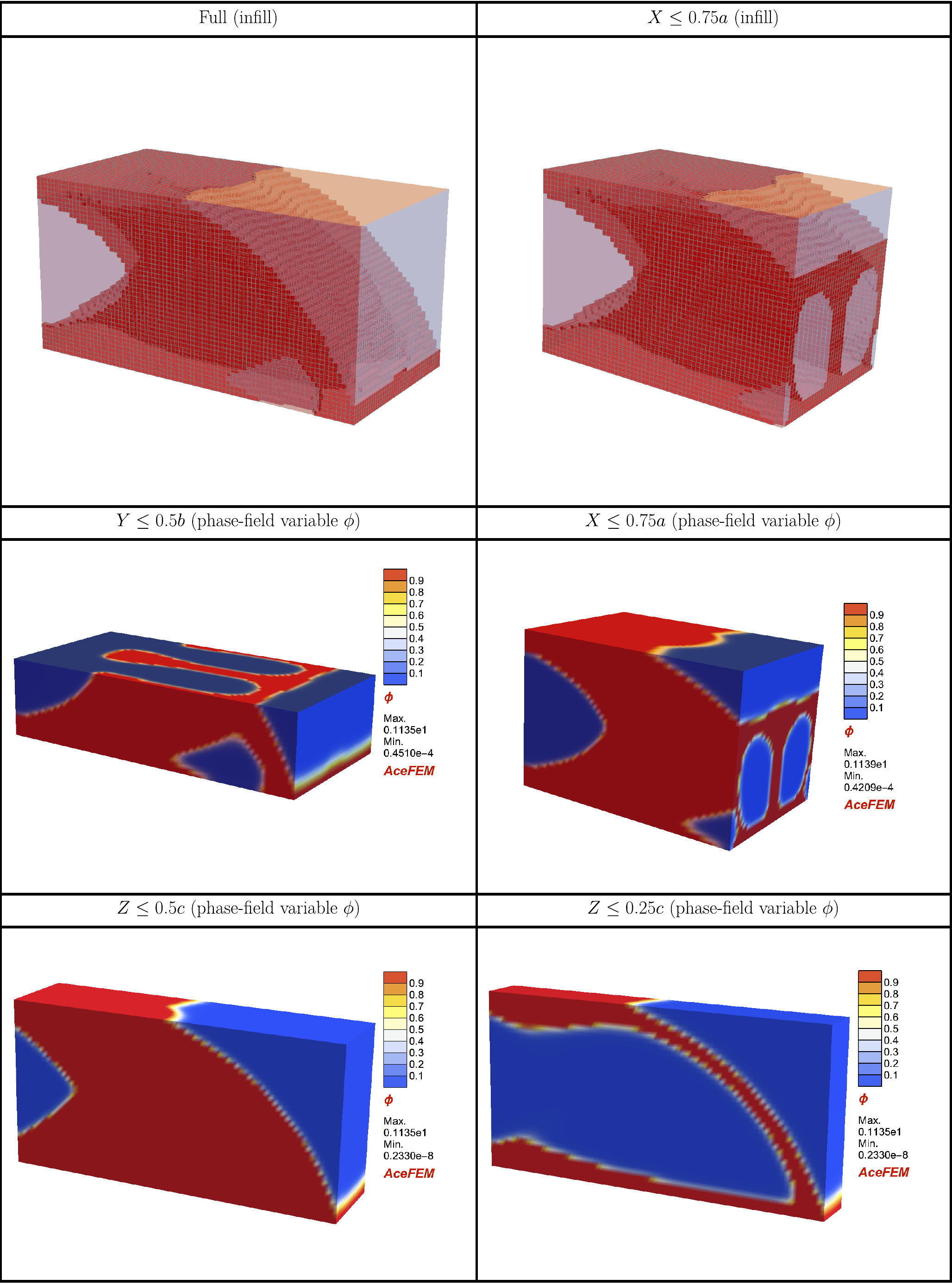}
\caption{Topology optimization of the 3D cantilever.}
\label{fig:3Dclamped}
\end{figure}

\begin{figure}[h!]
\centering
\includegraphics[width=0.75\textwidth]{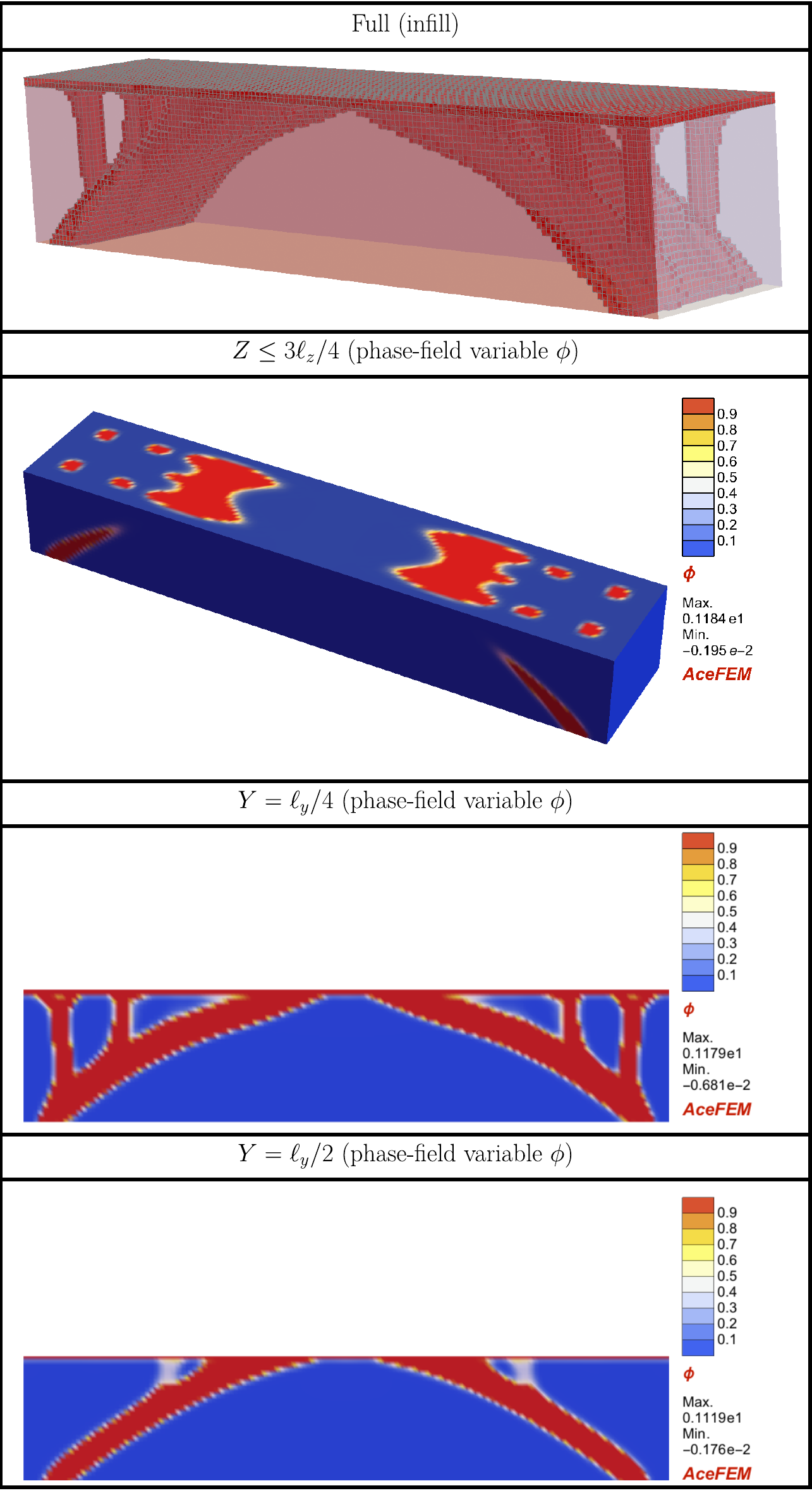}
\caption{Topology optimization of the 3D bridge.}
\label{fig:3Dbridge}
\end{figure}

The final obtained designs are shown in Fig. \ref{fig:3Dclamped} for the cantilever, and in Fig. \ref{fig:3Dbridge} for the bridge. Three-dimensional patterns clearly arise, with regular external and internal surfaces obtained without the need of post-processing filtering techniques thanks to the adopted phase-field rationale.

%% file: opt_conclusion_r1_v12.tex
\section{Comparison between investigated formulations and conclusions}

The present work addresses a rigorous study on basic properties of structural topology optimization problems with phase-field. In particular,  we performed the following five tasks:
\begin{enumerate}
\item Proposed mixed Hu-Washizu  variational formulations for the topology optimization problem with phase-field to directly impose equilibrium, constitutive, and compatibility equations in the formulation;
\item Investigated two different topology optimization principles, e.g.,: i) a formulation  imposing  \emph{a priori} the amount of material to be distributed within the design domain ({\em formulation with volume constraint}); ii) a formulation based on a minimization of material to be distributed, given that a {\it cost} (i.e., a penalty parameter) is assigned to the material ({\em formulation with volume minimization});
\item Introduced a Simultaneous Analysis and Design (SAND) monolithic solution strategy, thanks to the Hu-Washizu functional rationale and based on an Allen-Cahn scheme, where the phase-field variable evolves under the respect of mechanical equilibrium at each computational incremental step;
\item \UUU Analyzed \OOO both numerical convergence behaviours and obtained final designs, based on comparative analyses between simulation strategies (monolithic SAND vs. staggered NAND) and between topology optimization principles (volume constraint vs. minimization);
\item Analyzed the performance of the proposed variational formulation based on volume minimization also on three-dimensional case studies.
\end{enumerate}
On the basis of theoretical considerations, general guidelines are traced for properly setting the value of simulation parameters characterizing the phase-field evolution (i.e., $\tau_{\phi}$, $\kappa_{\phi}$ and $\gammaphi$). Following this careful setting, the implemented finite element formulations and solution algorithm are robust, allowing to conduct a wide campaign of parametric simulations without encountering any numerical convergence issue. 
In addition, the volume penalty parameter $\kappa_v$ in the volume minimization functional has been related to target values of volume fraction.

From the obtained results, we may conclude the following:

\begin{itemize}

\item A monolithic SAND solution strategy  significantly outperforms a staggered NAND solution strategy (see Fig. \ref{fig:conv_stagg_all});

\item The volume minimization formulation $\cL^{vm}$ shows an Allen-Cahn convergence behavior with better properties than the volume constraint formulation $\cL^{vc}$, since the decrease of the Allen-Cahn error measure $\mathcal{E}_{AC}$ is highly oscillatory for $\cL^{vc}$ but not for $\cL^{vm}$ (see Fig. \ref{fig:conv_all});

\item Comparing designs with the same final volume fraction  $v_{sol}$, the two formulations lead to practically coincident solutions in terms of compliance $\mathcal{C}_{sol}$ (or structural displacement $U_{sol}$) and interface perimeter $\mathcal{P}_{sol}$ (see Fig. \ref{fig:mass_2}), although final designs are slightly different (see Fig. \ref{fig:mass_maps});

\item The final design is practically independent from the applied load with the volume constraint functional $\cL^{vc}$, while is highly affected with the volume minimization one $\cL^{vm}$ (see Figs. \ref{fig:load} and \ref{fig:load_maps}),
which is seen as a significative advantage of this latter formulation;

\item The number of Newton-Raphson iterations required for solving the problem with the the volume minimization functional $\cL^{vm}$ is in most cases significantly lower (less than $50\%$) than the ones employed with the volume constraint functional $\cL^{vc}$;

\item No convergence issues are encountered in 3D applications (see Figs. \ref{fig:3Dclamped} and \ref{fig:3Dbridge}), allowing to obtain (thanks to the phase-field rationale) complex but regular patterns without the need of any post-processing filtering technique (as required by density-based approaches) which might lead to undesirable (artificial) effects \citep{Zegard2016}.

%\item converged values of compliance $\mathcal{C}_{sol}$ and interface perimeter $\mathcal{P}_{sol}$ have been obtained for  the volume minimization formulation $\cL^{vm}$ with higher tolerance (i.e., higher values of $c_{AC}^{conv}$) than for the volume constraint one $\cL^{vc}$, allowing to reduce the number of required Newton-Raphson iterations up to 60\% (cf., Figs. \ref{fig:CT1conv_all} and \ref{fig:CT0conv});

%\item the constraint $\phi \in [0,1]$ is more effectively enforced with the volume minimization formulation $\cL^{vm}$ than with the volume constraint one $\cL^{vc}$ (cf., Figs. \ref{fig:CT1mass_kb} and \ref{fig:CT0mass_kb}).

%allows to significantly reduce the number of Newton-Raphson iterations with respect to the volume constrained one in the range of final volume fractions herein investigated, and especially in the most challenging cases where $v_{sol} \in [0.2,0.5]$;

%\item  for the volume minimization formulation, it is possible to identify a relationship between the thickness penalization parameter $\gamma_{\phi}$ and the volume stiffness parameter $k_v$ (in turn related to the final volume fraction $v_{sol}$) which ensures an optimal stress distribution in the filled domain. Accordingly, the volume-minimization formulation seems to be favorable because of a higher robustness of final stress distributions, especially for low volume fractions. 

\end{itemize}

In conclusion, the obtained results highlight fundamental aspects of structural topology optimization problems with phase-field, allowing to optimize the solution strategy and to trace guidelines for settings simulation parameters beforehand. %, whose impact is often underestimated in the literature. On the other hand, 
Referring in particular to the proposed volume minimization principle, we believe that these results can be a starting point for more advanced developments of phase-field topology optimization that considers loading uncertainties \citep{Dunning2011} or multi-target strategies, e.g., controlling both geometry and compliance \citep{Stromberg2010}, both geometry and stresses \citep{burger_2006}, the structure life-cycle cost \citep{Sarma2002}, or manufacturing costs \citep{Liu2019}. 

\begin{acknowledgements}

M. Marino acknowledges partial support through the Rita Levi Montalcini Program for Young Researchers (Programma per Giovani Ricercatori $-$ anno 2017 Rita Levi Montalcini, Ministry of Education, University and Research, Italy).
F. Auricchio acknowledges partial support from  
 the Italian Minister of University and Research through the project A BRIDGE TO THE FUTURE: Computational methods, innovative applications, experimental validations of new materials and technologies (No. 2017L7X3CS) within the PRIN 2017 program and from Regione Lombardia through the project "MADE4LO - Metal ADditivE for LOmbardy" (No. 240963) within the POR FESR 2014-2020 program.
A. Reali acknowledges partial support from  
 the Italian Minister of University and Research through the project XFAST-SIMS (No. 20173C478N) within the PRIN 2017 program.
U. Stefanelli  acknowledges partial support through the FWF projects F\,65, I\,2375, and  P\,27052 and by the Vienna Science and Technology Fund (WWTF) through Project MA14-009.  
\end{acknowledgements}

%% file: opt_appendix_r1_v12.tex
\renewcommand\thefigure{\thesection.\arabic{figure}}   
\renewcommand\thetable{\thesection.\arabic{table}}   

\setcounter{figure}{0}    
\setcounter{subfigure}{0}
\setcounter{equation}{0}

\section{Interpolation of the stress and strain fields} \label{app:interpolation}

The interpolation matrices of the stress and strain fields employed for 2D and 3D applications are detailed in what follows. Choices follow arguments presented in \cite{Weissman96}, \cite{Cao2002} and  \cite{Djoko2006}. In detail, minimum distributions required for stability are included as a common basis between the assumed stress and strain fields. In addition, the assumed strain field is enriched with respect to the stress one with strain modes not already contained in the minimum one. A minimal strain enrichment, \emph{a priori} orthogonal to the assumed stress field, is considered. Interesting extensions of the present work would be to investigate more refined enrichment strategies, accompanied by a generalization of the introduced mixed formulations, in order to deal with material constrained behaviours (e.g., incompressibility or inextensibility). 

\subsection{2D applications}

The stress field $\ffs$  is approximated in a discontinuous element-by-element form, introducing five element unknowns collected in  vector $\hat{\ffs}=(\hat{\sigma}_1,\ldots,\hat{\sigma}_5)$ and the following definition of $\bbN_{\ffs}(\ffxi)$:
\begin{equation}
\text{vec}(\ffs_{\ffxi}) = \left(
\begin{array}{c}
\sigma_{\xi \xi} \\
\sigma_{\eta\eta} \\
\sigma_{\xi\eta}
\end{array}
\right) =
\bbN_{\ffs}(\ffxi) \hat{\ffs}
\htext{3mm}{with}
\bbN_{\ffs}(\ffxi) =
\left[
\begin{array}{ccccccc}
1 & \eta & 0 & 0  & 0 \\
0 & 0     & 1 & \xi & 0\\
0 & 0     & 0 &   0 & 1
\end{array}
\right]\, .
\end{equation}

The strain field $\ffe$  is approximated in a discontinuous element-by-element form, introducing seven element unknowns collected in vector $\hat{\ffe}=(\hat{\epsilon}_1,\ldots,\hat{\epsilon}_7)$ and the following definition of $\bbN_{\ffe}(\ffxi) $
\begin{equation}
\text{vec}(\ffe_{\ffxi}) = \left(
\begin{array}{c}
\epsilon_{\xi \xi} \\
\epsilon_{\eta\eta} \\
2 \epsilon_{\xi\eta}
\end{array}
\right) =
\bbN_{\ffe}(\ffxi) \hat{\ffe}
\htext{3mm}{with}
\bbN_{\ffe}(\ffxi) =
\left[
\begin{array}{ccccccc}
1 & \xi & \eta& 0 & 0 &  0 & 0 \\
0 & 0 & 0     & 1 & \xi & \eta & 0\\
0 & 0 & 0 & 0 & 0 & 0 & 1
\end{array}
\right]\, .
\end{equation}

\subsection{3D applications}

The stress field $\ffs$  is approximated in a discontinuous element-by-element form, introducing eighteen element unknowns collected in  vector $\hat{\ffs}=(\hat{\sigma}_1,\ldots,\hat{\sigma}_{18})$ and the following definition of $\bbN_{\ffs}(\ffxi)$:
\begin{equation}
\text{vec}(\ffs_{\ffxi}) = \left(
\begin{array}{c}
\sigma_{\xi \xi} \\
\sigma_{\eta\eta} \\
\sigma_{\zeta \zeta}\\
\sigma_{\eta\zeta}\\
\sigma_{\xi\zeta} \\
\sigma_{\xi\eta}
\end{array}
\right) =
\bbN_{\ffs}(\ffxi) \hat{\ffs}
\htext{3mm}{with}
\bbN_{\ffs}(\ffxi) =
\left[
\begin{array}{cccccccccccccccccc}
1 & \eta & \zeta & \eta \zeta & 0 & 0  & 0 & 0 & 0  & 0 & 0 & 0  & 0 & 0 & 0 & 0  & 0 & 0    \\
0 & 0 & 0 & 0 & 1 & \xi & \zeta & \xi \zeta & 0 & 0 & 0  & 0 & 0 & 0  & 0 & 0 & 0  & 0 \\
0 & 0     & 0 &   0 & 0 & 0 & 0 & 0 & 1 & \xi & \eta & \xi\eta & 0 & 0  & 0 & 0 & 0  & 0 \\
0 & 0     & 0 &   0 & 0 & 0 & 0 & 0 & 0 & 0 & 0 & 0 &1 & \xi  &  0 & 0 & 0  & 0 \\
0 & 0     & 0 &   0 & 0 & 0 & 0 & 0 & 0 & 0 & 0 & 0 &  0 & 0 &1 & \eta & 0  & 0 \\
0 & 0     & 0 &   0 & 0 & 0 & 0 & 0 & 0 & 0 & 0 & 0 &  0 & 0  & 0  & 0 &1 & \zeta
\end{array}
\right]\, .
\end{equation}

The strain field $\ffe$  is approximated in a discontinuous element-by-element form, introducing twenty-one element unknowns collected in vector $\hat{\ffe}=(\hat{\epsilon}_1,\ldots,\hat{\epsilon}_{21})$ and the following definition of $\bbN_{\ffe}(\ffxi) $
\begin{equation}
\text{vec}(\ffe_{\ffxi}) = \left(
\begin{array}{c}
\epsilon_{\xi \xi} \\
\epsilon_{\eta\eta} \\
\epsilon_{\zeta \zeta}\\
2\epsilon_{\eta\zeta}\\
2\epsilon_{\xi\zeta} \\
2\epsilon_{\xi\eta}
\end{array}
\right) =
\bbN_{\ffs}(\ffxi) \hat{\ffs}
\htext{3mm}{with}
\bbN_{\ffe}(\ffxi) =
\left[
\begin{array}{ccccccccccccccccccccc}
1 & \xi & \eta & \zeta & \eta \zeta & 0 & 0 & 0 & 0 & 0 & 0 & 0  & 0 & 0 & 0  & 0 & 0 & 0 & 0  & 0 & 0    \\
0 & 0 & 0 & 0 & 0 & 1 & \xi & \eta & \zeta & \xi \zeta & 0 & 0& 0 & 0  & 0 & 0 & 0  & 0 & 0 & 0  & 0 \\
0 & 0  & 0   & 0 &   0 & 0 & 0 & 0 & 0 & 0 & 1 & \xi & \eta & \zeta & \xi\eta & 0 & 0  & 0 & 0 & 0  & 0 \\
0 & 0   & 0 & 0 & 0  & 0 &   0 & 0 & 0 & 0 & 0 & 0 & 0 & 0 & 0 &1 & \xi  &  0 & 0 & 0  & 0 \\
0 & 0  & 0 & 0 & 0   & 0 &   0 & 0 & 0 & 0 & 0 & 0 & 0 & 0 & 0 &  0 & 0 &1 & \eta & 0  & 0 \\
0 & 0   & 0 & 0 & 0  & 0 &   0 & 0 & 0 & 0 & 0 & 0 & 0 & 0 & 0 &  0 & 0  & 0  & 0 &1 & \zeta
\end{array}
\right]\, .
\end{equation}

\setcounter{figure}{0}    
\setcounter{subfigure}{0}
\setcounter{equation}{0}

%% file: Articolo_r1_v12.bbl
\begin{thebibliography}{29}
\expandafter\ifx\csname natexlab\endcsname\relax\def\natexlab#1{#1}\fi
\expandafter\ifx\csname url\endcsname\relax
  \def\url#1{\texttt{#1}}\fi
\expandafter\ifx\csname urlprefix\endcsname\relax\def\urlprefix{URL }\fi

\bibitem[{Allaire et~al.(2004)Allaire, Jouve, and Maillot}]{allaire_2004}
Allaire, G., Jouve, F., Maillot, H., 2004. Topology optimization and optimal
  shape design using homogenization. Struct. Multidisc. Optim. 28, 87--98.

\bibitem[{Auricchio et~al.(2019)Auricchio, Bonetti, Carraturo, H\"{o}mberg,
  Reali, and E.}]{auricchio_19}
Auricchio, F., Bonetti, E., Carraturo, M., H\"{o}mberg, D., Reali, A., E., R.,
  2019. A phase-field based on graded-material topology optimization with
  stress constraint. arXiv:1907.06355, Math. Models Meth. Appl. Sci., to
  appear.

\bibitem[{Bends{\o}e(1983)}]{bendsoe_83}
Bends{\o}e, M.~P., 1983. On obtaining a solution to optimization problems for
  solid, elastic plates by restriction of the design space. J. Struct. Mech.
  11~(4), 501--521.

\bibitem[{Bends{\o}e and Sigmund(1999)}]{bendsoe_99}
Bends{\o}e, M.~P., Sigmund, O., 1999. Material interpolation schemes in
  topology optimization. Archive of Applied Mechanics 6 65, 635--654.

\bibitem[{Blank et~al.(2014{\natexlab{a}})Blank, Farshbaf-Shaker, Garcke,
  Rupprecht, and Styles}]{blank_multi-material_2014}
Blank, L., Farshbaf-Shaker, M., Garcke, H., Rupprecht, C., Styles, V.,
  2014{\natexlab{a}}. Multi-material {Phase} {Field} {Approach} to {Structural}
  {Topology} {Optimization}. In: Leugering, G., Benner, P., Engell, S.,
  Griewank, S., Harbrecht, H., Hinze, M., Rannacher, R., Ulbrich, S. (Eds.),
  Trends in {PDE} {Constrained} {Optimization}. Vol. 165. Springer
  International Publishing, Cham, pp. 231--246.

\bibitem[{Blank et~al.(2014{\natexlab{b}})Blank, Garcke, Farshbaf-Shaker, and
  V.}]{blank_14}
Blank, L., Garcke, H., Farshbaf-Shaker, M., V., S., 2014{\natexlab{b}}.
  Relating phase field and sharp interface approaches to structural topology
  optimization. ESAIM Control Optim. Calc. Var. 20, 1025--1058.

\bibitem[{Bourdin and Chambolle(2003)}]{bourdin_2003}
Bourdin, B., Chambolle, A., 2003. Design-dependent loads in topology
  optimization. ESAIM Contr. Optim. Calc. Var. 9, 19--48.

\bibitem[{Burger(2003)}]{burger_2003}
Burger, M., 2003. A framework for the construction of level set methods for
  shape optimization and reconstruction. Interfaces Free Bound. 5, 301--332.

\bibitem[{Burger and Stainko(2006)}]{burger_2006}
Burger, M., Stainko, R., 2006. Phase-field relaxation of topology optimization
  with local stress constraints. SIAM J. Control Optim. 45~(4), 1447--1466.

\bibitem[{Cao et~al.(2002)Cao, Hu, Lu, Fukunaga, and Yao}]{Cao2002}
Cao, Y.~P., Hu, N., Lu, J., Fukunaga, H., Yao, Z.~H., 2002. A {3D} brick
  element based on {H}u-{W}ashizu variational principle for mesh distortion.
  Internat. J. Numer. Methods Engrg. 53~(11), 2529--2548.

\bibitem[{Carraturo et~al.(2019)Carraturo, Rocca, Bonetti, H\"{o}mberg, Reali,
  and F.}]{carraturo_19}
Carraturo, M., Rocca, E., Bonetti, E., H\"{o}mberg, D., Reali, A., F., A.,
  2019. Graded-material {D}esign based on {P}hase-field and {T}opology
  {O}ptimization. Computational Mechanics, DOI: 10.1007/s00466--019--01736--w.

\bibitem[{Deaton and Grandhi(2014)}]{Deaton2014}
Deaton, J.~D., Grandhi, R.~V., 2014. A survey of structural and
  multidisciplinary continuum topology optimization: post 2000. Structural and
  Multidisciplinary Optimization 49~(1), 1--38.

\bibitem[{Ded{\`e} et~al.(2012)Ded{\`e}, Borden, and
  Hughes}]{dede_isogeometric_2012}
Ded{\`e}, L., Borden, M.~J., Hughes, T.~J., 2012. Isogeometric analysis for
  topology optimization with a phase field model. Archives of Computational
  Methods in Engineering 19~(3), 427--465.

\bibitem[{Djoko et~al.(2006)Djoko, Lamichhane, Reddy, and Wohlmuth}]{Djoko2006}
Djoko, J., Lamichhane, B., Reddy, B., Wohlmuth, B., 2006. Conditions for
  equivalence between the {H}u-{W}ashizu and related formulations, and
  computational behavior in the incompressible limit. Computer Methods in
  Applied Mechanics and Engineering 195~(33), 4161 -- 4178.

\bibitem[{Dunning et~al.(2011)Dunning, Kim, and Mullineux}]{Dunning2011}
Dunning, P.~D., Kim, H.~A., Mullineux, G., 2011. Introducing loading
  uncertainty in topology optimization. AIAA Journal 49~(4), 760--768.

\bibitem[{Korelc and Wriggers(2016)}]{Korelc2016}
Korelc, J., Wriggers, P., 2016. Automation of Finite Element Methods. Springer
  International Publishing.

\bibitem[{Liu et~al.(2019)Liu, Chen, Liang, and To}]{Liu2019}
Liu, J., Chen, Q., Liang, X., To, A.~C., 2019. Manufacturing cost constrained
  topology optimization for additive manufacturing. Frontiers of Mechanical
  Engineering 14~(2), 213--221.

\bibitem[{Modica(1987)}]{modica_87}
Modica, L., 1987. The gradient theory of phase transitions and the minimal
  lnterface criterion. Arch. Rat. Mech. Anal.. 98, 123--142.

\bibitem[{Osher and Santosa(2001)}]{osher_santosa_01}
Osher, S., Santosa, F., 2001. Level set methods for optimization problems
  involving geometry and constraints i. frequencies of a two-density
  inhomogeneous drum. J. Comput. Phys. 171, 272--288.

\bibitem[{Penzler et~al.(2012)Penzler, Rumpf, and Wirth}]{penzler_2012}
Penzler, P., Rumpf, M., Wirth, B., 2012. A phase-field model and minimal
  compliance shape optimization in nonlinear elasticity. ESAIM Control Optim.
  Calc. Var. 2012, 229--258.

\bibitem[{Pian and Sumihara(1984)}]{PianSumihara1984}
Pian, T. H.~H., Sumihara, K., 1984. Rational approach for assumed stress finite
  elements. International Journal for Numerical Methods in Engineering 20~(9),
  1685--1695.

\bibitem[{Sarma and Adeli(2002)}]{Sarma2002}
Sarma, K.~C., Adeli, H., 2002. Life-cycle cost optimization of steel
  structures. International Journal for Numerical Methods in Engineering
  55~(12), 1451--1462.

\bibitem[{Sigmund and Petersson(1998)}]{sigmund_98}
Sigmund, O., Petersson, J., 1998. Numerical instabilities in topology
  optimization: A survey on procedures dealing with cheackboards,
  mesh-dependencies and local minima. Structural Optimization 16, 68--75.

\bibitem[{Str{\"o}mberg(2010)}]{Stromberg2010}
Str{\"o}mberg, N., 2010. Topology optimization of structures with manufacturing
  and unilateral contact constraints by minimizing an adjustable
  compliance--volume product. Structural and Multidisciplinary Optimization
  42~(3), 341--350.

\bibitem[{Suzuki and Kikuchi(1991)}]{suzuki_91}
Suzuki, K., Kikuchi, N., 1991. A homogenization method for shape and topology
  optimization. Comput. Methods Appl. Mech. Engrg. 93~(3), 291--318.

\bibitem[{Takezawa et~al.(2010)Takezawa, Nishiwaki, and
  Kitamura}]{takezawa_2010}
Takezawa, A., Nishiwaki, S., Kitamura, M., 2010. Shape and topology
  optimization based on the phase field method and sensitivity analysis. J.
  Comp. Phys. 229~(7), 2697--2718.

\bibitem[{Weisman(1996)}]{Weissman96}
Weisman, S.~L., 1996. High-accuracy low-order {T}hree-dimensional {B}rick
  {E}lements. International Journal for Numerical Methods in Engineering
  39~(14), 2337--2361.

\bibitem[{Zegard and Paulino(2016)}]{Zegard2016}
Zegard, T., Paulino, G.~H., 2016. Bridging topology optimization and additive
  manufacturing. Structural Multidisciplinary Optimization 53, 175--192.

\bibitem[{Zhou and Rozvany(1991)}]{zhou_91}
Zhou, M., Rozvany, G. I.~N., 1991. The coc algorithm, part ii: Topological
  geometry and generalized shape optimization. Comp. Meth. Appl. Mech. Engng.
  89, 197--224.

\end{thebibliography}
